\documentstyle[a4,amssymb,12pt]{amsart}

\newtheorem{theorem}{Theorem}[section]

\newtheorem{lemma}[theorem]{Lemma}
\newtheorem{proposition}[theorem]{Proposition}
\newtheorem{definition}{Definition}[section]
\newtheorem{remark}{Remark}[section]
\numberwithin{equation}{section}

\begin{document}
\title[q-deformed universal DS-reduction]
{Q-pseudodifference Drinfeld-Sokolov reduction for 
algebra of complex size matrices.}
\author{A.L.Pirozerski and M.A.Semenov-Tian-Shansky}
\maketitle

\begin{center}
Universit\'e de Bourgogne, Dijon, France \\[0pt]
and Steklov Institute of Mathematics, St.Petersburg, Russia
\end{center}

\begin{abstract}
The q-deformed version of the Drinfeld-Sokolov reduction is extended
to the case of the algebra of 'complex size matrices'; this 
construction generalizes earlier results of B.Khesin and F.Malikov
on universal DS reduction and follows the pattern of recent studies
of q-deformed DS reduction which were started by E.Frenkel, N.Reshetikhin
and one of the authors.
\end{abstract}
\bigskip
\section{Introduction}

It is well-known that the space ${\cal M}_{n}$ of scalar $n$-th order
differential operators has a remarkable quadratic Poisson structure, called
the (second) Adler-Gelfand-Dickey bracket \cite{a, gd}. This structure
admits several different realizations. The first one, known as the
Drinfeld-Sokolov reduction \cite{ds}, shows that the Gelfand-Dickey bracket
can be obtained via Hamiltonian reduction from a linear Poisson bracket on
the space of matrix first order differential operators which is considered
as the dual space of the affine algebra $\widehat{{\frak sl}}_{n}.$ This
construction has a natural generalization to the case of an arbitrary
semisimple Lie algebra ${\frak g}$, the corresponding Poisson algebra ${\cal %
W}({\frak g})$ of functionals on reduced space is called the classical $%
{\cal W}$-algebra (associated with ${\frak g}$).

The second realization is based on the study of the Lie group of integral
operators (more precisely, of some its extension, see \cite{KhZ}). This
group comes equipped with the natural Sklyanin bracket which endows it with
the structure of a Poisson-Lie group, and the second Gelfand-Dickey bracket
is identified with the restriction of this bracket to the subvariety ${\cal {%
M}}_n,$ which is Poisson in this case.

This construction leads to a natural generalization of the Gelfand-Dickey
bracket for the space of pseudodifferential symbols of any complex degree $%
\lambda $ 
\begin{equation*}
G_{\lambda }=\left\{ \partial ^{\lambda }+u_{1}\left( x\right) \partial
^{\lambda -1}+\ldots \right\} .
\end{equation*}
At the same time, DS-reduction \cite{ds} is defined only for $\lambda =n,$ $%
n\in {\Bbb N},$ and only for the Poisson subspace ${\cal M}_{n}\subset
G_{n}. $ B.Khesin and F.Malikov in \cite{KhM} proposed a counterpart of the
DS-reduction which applies to all pseudodifferential symbols of complex
degree. To describe their construction (called the universal DS-reduction)
let us recall the definition of the algebras ${\frak gl}_{\lambda }$ of
complex size matrices. The definition of ''complex size matrices'' was
proposed by B.Feigin, who used them to compute the cohomology of the Lie
algebra of differential operators, see \cite{feig}.

Consider the universal enveloping algebra $U\left( {\frak {sl}}_{2} \right) $%
; it has the well-known Casimir element $C=ef+fe+\frac{1}{2}h^{2},$ where $%
e,f,h$ --- standard basis of ${\frak {sl}}_{2}.$ By definition, 
\begin{equation*}
{\frak gl}_{\lambda }=U\left( {\frak {sl}}_{2}\right) /\left\{ C=\frac{1}{2}%
(\lambda -1)(\lambda +1)\right\} .
\end{equation*}
For positive integer $\lambda =n$ the algebra ${\frak {gl}}_{\lambda }$
contains a large ideal, and the quotient algebra is isomorphic to ${\frak {gl%
}}_{n},$ which explains the name ''algebra of complex size matrices''.

More precisely, to perform the universal DS-reduction we need not the
algebra ${\frak {gl}}_{\lambda }$ itself, but some its completion $\bar{%
{\frak {gl}}}_{\lambda }.$ B.Khesin and F.Malikov \cite{KhM} proved that the
Gelfand-Dickey bracket on $G_{\lambda }$ can be obtained by reduction from a
linear bracket on the dual of the affine Lie algebras corresponding to $\bar{%
{\frak {gl}}}_{\lambda }.$ (This was conjectured earlier by B.Feigin and
C.Roger).

The third realization of Gelfand-Dickey brackets is based on the study of
the center of the quantized enveloping algebra $U_q\left( \widehat{{\frak sl}}%
_{n}\right) $ at the critical value of the central charge, see \cite{ef}. We
shall not consider it here, although it is this realization which was
generalized for the first time to q-difference setting and allowed to
construct the q-deformed ${\cal W}$-algebras (see \cite{FrRes}).

A q-difference version of the DS reduction was defined in \cite{FRS} for the 
${\frak sl}_{2}$ case and generalized to the case of arbitrary semisimple
algebra in \cite{SemSev}. The consistency conditions for this reduction lead
to a new class of elliptic r-matrices which are fixed by these conditions 
in an essentially unique way (one for each semisimple Lie algebra).

The q-difference counterpart of the second construction of the
Gelfand-Dickey brackets was proposed in \cite{PirSem}. It is based on the
study of the group of q-pseudodifference symbols of arbitrary complex
degrees 
\begin{equation*}
\widehat{{\Bbb G}}_{-}=\bigcup_{\lambda \in {\Bbb C}}\widehat{{\Bbb G}}%
_{\lambda },\quad \widehat{{\Bbb G}}_{\lambda }=\left\{ L=D^{\lambda
}+u_{1}\left( z\right) D^{\lambda -1}+\cdots \right\} ,
\end{equation*}
where $D$ is the dilation operator, $(Df)(z)=f(qz).$ It was shown that in a
natural class of r-matrix Poisson brackets on $\widehat{{\Bbb G}}_{\lambda }$
there exists a unique one with respect to which formal spectral invariants $%
H_{n}(L)=\frac{\lambda }{n}{\rm Tr\,}L^{\frac{n}{\lambda }}$ are in
involution. (These spectral invariants give rise to the generalized
q-deformed KdV hierarchy $\frac{dL}{dt}=\left[ L_{\left( +\right) }^{\frac{m%
}{\lambda }},L\right] $ described earlier in \cite{Frenkel}.) This bracket
was called the generalized q-deformed Gelfand-Dickey structure. Similarly to
the differential case, for positive integer $\lambda =n$ the subspace of
q-difference operators of $n$-th order 
\begin{equation*}
{\Bbb M}_{n}=\left\{ L=D^{n}+u_{1}\left( z\right) D^{n-1}+\cdots
+u_{0}\left( z\right) \right\}
\end{equation*}
is Poisson. The Poisson algebra of the functionals on ${\Bbb M}_{n}$
coincides (up to the constraint $u_{0}\left( z\right) =1$) with the
q-deformed ${\cal W}$-algebra ${\cal W}_{q}\left( {\frak sl}_{n}\right) $,
constructed in \cite{FrRes} and \cite{FRS, SemSev}.

In the present paper we generalize to the q-difference setting the procedure
of the universal DS-reduction \cite{KhM}. We define an algebra ${\frak gl}%
_{q}$ consisting of ${\frak gl}_{\infty }$-matrices of special form whose
matrix elements are holomorphic functions of a complex variable $t.$ (This
algebra is an extension of the algebra ${\frak gl},$ which was also
constructed in \cite{KhM}). There exists a natural evaluation map $%
i_{\lambda }:{\frak {gl}}_{q}\mapsto {\frak {gl}}_{\infty },$ which assigns
to a matrix-function $A\in {\frak gl}_{q}$ its value at an arbitrary fixed
point $t=\lambda $: $A\mapsto A(\lambda ).$ The image of this evaluation map
is a subalgebra in ${\frak gl}_{\infty }$ which will be called (extended)
algebra of complex size matrices (more precisely, of size $\lambda \times
\lambda $. In fact, we have a family of such algebras parametrized by a
complex number $\lambda $). Due to the infinite dimension of ${\frak gl}%
_{q}^{\lambda },$ the definition of the corresponding loop algebras $L{\frak %
gl}_{q}^{\lambda }$ involves some peculiarities.

Like in the ${\frak sl}_{n}$-case \cite{FRS, SemSev}, the reduction
procedure consists of two steps: first, we impose constraints fixing a
submanifold ${{\Bbb {Y}}}_{q}^{\lambda }\subset L{\frak gl}_{q}^{\lambda }$
which is preserved by the q-deformed gauge action of the upper-triangular
group $L{{\frak {N}}}_{+},$ and then we take the quotient over this group.
We prove that the quotient ${{\Bbb {Y}}}_{q}^{\lambda }/L{{\frak {N}}}_{+}$
can be identified with the space $\widehat{{\Bbb G}}_{\lambda }$ of
q-pseudodifference symbols of degree $\lambda .$

Using a method similar to the one of \cite{FRS, SemSev}, we describe
explicitly all r-matrix Poisson brackets on $L{\frak gl}_{q}^{\lambda }$ (in
a wide natural class) which admit the q-deformed universal DS-reduction. At
this point we encounter a new phenomenon. Recall, that in the ${\frak sl}%
_{n} $-case the underlying classical r-matrix was related with the
decomposition of the algebra $L{\frak sl}_{n}$ into the sum of the
subalgebras of upper-triangular, lower-triangular and diagonal matrices. Its
diagonal part $\hat{r}^{0}$ was given by the Cayley transformation of the
operator $D\tau _{n},$ the operator $\tau _{n}$ acting by {\em cyclic
permutation} of matrix elements. Obviously, in the ${\frak gl}_{q}^{\lambda
} $-case there exists no analogue of $\tau _{n}.$ We shall see that it must
be replaced by {\em the shift operator} $\hat{s}$ whose properties are quite
different. This causes some difficulties in the definition of the diagonal
part of the r-matrix, which requires a regularization introducing some free
parameters into the admissible r-matrix; the uniqueness is restored,
however, if we demand the formal spectral invariants to be in involution,
and the resulting quotient Poisson structure coincides precisely with the
generalized q-deformed Gelfand-Dickey structure defined in \cite{PirSem}.

This article has the following structure. In section 2 we recall in more
details the finite-dimensional DS-reduction \cite{FRS, SemSev}, as well as
some results of \cite{PirSem} and generalize them to the ${\frak gl}_{n}$%
-case. Unlike the ${\frak sl}_{n}$-case, we obtain a family of Poisson
structures admitting reduction; however, only one of them gives rise to a
quotient bracket satisfying the involutivity condition. This result is a
finite-dimensional analogue of the uniqueness theorem for ${\frak {gl}}%
_{q}^{\lambda }$ and will be used to prove the latter.

In section 3 we construct algebras ${\frak gl}_{q},$ ${\frak gl}%
_{q}^{\lambda }$ and their loop algebras. Section 4 is devoted to the
cross-section theorem which gives a model of the quotient ${{\Bbb {Y}}}%
_{q}^{\lambda }/L{{\frak {N}}}_{+}$. In section 5 we describe explicitly all
r-matrix Poisson brackets on $L{\frak {gl}}_{q}^{\lambda },$ which admit the
q-deformed universal DS-reduction. Section 6 is devoted to the uniqueness
theorem.

Throughout the article we shall use the following notation. We fix a complex
number $q,$ $|q|<1.$ Let $\hat{h}$ be the dilation operator, 
\begin{equation}
\hat{h}a(z)=a(qz),\quad a\in {{\Bbb {C}}}\left( \left( z^{-1}\right) \right)
,\quad q\in {{\Bbb {C}}},\quad \left| q\right| <1.  \label{e010}
\end{equation}
We shall denote the same dilation operator by $D$ when it is considered as a
generator of the algebra of q-pseudodifference operators (see below).

We fix the branch of $\ln w,$ $w\in {{\Bbb {C}}},$ by 
\begin{equation*}
-\pi <\arg w<\pi ,\quad \ln 1=0,
\end{equation*}
and put $q^{w}\equiv \exp \left( w\ln q\right) .$ Arbitrary complex degrees
of $\hat{h}$ are defined by 
\begin{equation}
^{h^{w}}a(z)\equiv \left( \hat{h}^{w}a\right) (z)=a\left( q^{w}z\right)
,\quad \forall w\in {{\Bbb {C}},\quad }a\in {{\Bbb {C}}}\left( \left(
z^{-1}\right) \right) .  \label{e011}
\end{equation}
For $a\in {{\Bbb {C}}}\left( \left( z^{-1}\right) \right) ,\quad
a=\sum_{i}a_{i}z^{i},$ we put 
\begin{equation}
\int a(z)dz/z=\mathop{Res}a=a_{0};  \label{e012}
\end{equation}
clearly, this formal integral is dilation invariant, i.e., 
\begin{equation}
\int a(z)dz/z=\int a(qz)dz/z.  \label{e013}
\end{equation}
We introduce an $\hat{h}$-invariant inner product in ${{\Bbb {C}}}\left(
\left( z^{-1}\right) \right) $ by 
\begin{equation}
\left\langle a,b\right\rangle _{{{\Bbb {C}}}}=\int \frac{dz}{z}a(z)b(z).
\label{e014}
\end{equation}

Let ${{\frak {a}}}$ be a linear space with an inner product $\left\langle
\cdot ,\cdot \right\rangle .$ We shall denote by $\left\langle \left\langle
\cdot ,\cdot \right\rangle \right\rangle $ the following inner product on ${%
{\frak {a}}}\oplus {{\frak {a}}}$: 
\begin{equation}
\left\langle \left\langle \left( 
\begin{array}{l}
X_{1} \\ 
X_{2}
\end{array}
\right) \,,\left( 
\begin{array}{l}
Y_{1} \\ 
Y_{2}
\end{array}
\right) \right\rangle \right\rangle =\left\langle X_{1},Y_{1}\right\rangle
-\left\langle X_{2},Y_{2}\right\rangle .  \label{e015}
\end{equation}

\section{The q-deformed DS-reduction in ${\frak gl}_{n}$: overview of the
results \protect\cite{FRS, SemSev, PirSem} and some generalizations.}

\subsection{Reduction procedure and the choice of r-matrix.}

In this subsection we recall briefly the procedure of DS-reduction and
describe all the r-matrix Poisson brackets on ${\frak gl}_{n}$ (in a wide
natural class) which admit this reduction. Unlike the ${\frak sl}_{n}$-case
where this bracket is essentially unique (see \cite{FRS, SemSev}), there
exists a family of such brackets parametrized by a skew-symmetric operator
in ${{\Bbb {C}}}\left( \left( z^{-1}\right) \right) .$ This non-uniqueness
does not lead to a new kind of deformed Gelfand-Dickey brackets, since, as
we shall see in subsection 2.2, only one of the quotient Poisson structures
satisfies the involutivity condition.

We fix the following notation. Let ${{\frak {n}}}_{+}\left( n\right) ,$ ${%
{\frak {n}}}_{-}\left( n\right) ,$ ${{\frak {h}}}_{n}\subset {\frak gl}_{n}$
be the subalgebras of strictly upper triangular, strictly lower triangular
and diagonal matrices, respectively, let ${{\frak {b}}}_{\pm }\left(
n\right) ={{\frak {n}}}_{\pm }\left( n\right) \oplus {{\frak {h}}}_{n}$ be
the subalgebras of upper (lower) triangular matrices with arbitrary diagonal
elements; let ${\bf N}_{+}\left( n\right) \subset {\frak gl}_{n}$ be the
unipotent group corresponding to ${{\frak {n}}}_{+}\left( n\right) .$ We
shall denote by $L{\frak gl}_{n},$ $L{{\frak {n}}}_{\pm }\left( n\right) ,$ $%
L{{\frak {b}}}_{\pm }\left( n\right) ,$ $L{{\frak {h}}}_{n},$ $L{\bf N}%
_{+}\left( n\right) $ the corresponding loop algebras (group).

We introduce an invariant inner product on $L{\frak gl}_{n}$ by 
\begin{equation}
\left\langle A,B\right\rangle =\int \frac{dz}{z}{\rm \mathop{Tr}}A\left(
z\right) B\left( z\right) ,\quad A,B\in L{\frak gl}_{n}.  \label{e11}
\end{equation}

We shall denote by ${{\Bbb {M}}}_{n}$ the space of scalar $n$-th order
q-difference operators of the form 
\begin{equation}
L=D^{n}+u_{1}(z)D^{n-1}+\cdot \cdot \cdot +u_{n}(z),\quad u_{i}\in {{\Bbb {C}%
}}\left( \left( z^{-1}\right) \right) .  \label{e19}
\end{equation}

Now let us briefly recall the DS-reduction procedure. It is well known that
a scalar q-difference equation of order $n$ 
\begin{equation*}
L\psi _{0}=0,\quad L\in {{\Bbb {M}}}_{n},
\end{equation*}
is equivalent to a first order matrix equation 
\begin{equation*}
D\Psi ={{\Bbb {L}}}\Psi ,\quad \Psi =\left( 
\begin{array}{c}
\psi _{n-1} \\ 
\vdots \\ 
\psi _{0}
\end{array}
\right) ,
\end{equation*}
where the potential ${{\Bbb {L}}}\in L{\frak gl}_{n}$ has a special form.
The standard choice for ${{\Bbb {L}}}$ is given by a companion matrix: 
\begin{equation}
{{\Bbb {L}}}=\left( 
\begin{array}{cccc}
-u_{1} & \cdots & -u_{n-1} & -u_{n} \\ 
1 & \ddots & 0 & 0 \\ 
\vdots & \ddots & \ddots & \vdots \\ 
0 & \cdots & 1 & 0
\end{array}
\right) .  \label{e110}
\end{equation}
This choice is not unique; a linear change of variables 
\begin{equation*}
\Psi \mapsto \Psi ^{\prime }=S\Psi ,\quad S\in L{\bf N}_{+}\left( n\right)
\end{equation*}
induces a gauge transformation 
\begin{equation}
{{\Bbb {L}}} \longmapsto {{\Bbb {L}}}^{\prime }=\,^{h}S{{\Bbb {L}}}S^{-1}.
\label{e111}
\end{equation}
Let us denote by ${{\Bbb {Y}}}_{n}\subset L{\frak gl}_{n}$ the subvariety of
all matrices of the form ${{\Bbb {L}}}^{\prime }=\Lambda _{n}+A,$ $A\in L{%
{\frak {b}}}_{+}\left( n\right) ,$ where 
\begin{equation}
\Lambda _{n}=\left( 
\begin{array}{cccc}
0 & \cdots & 0 & 0 \\ 
1 & \ddots & 0 & 0 \\ 
\vdots & \ddots & \ddots & \vdots \\ 
0 & \cdots & 1 & 0
\end{array}
\right) .  \label{e112}
\end{equation}
It is easy to see that the gauge action (\ref{e111}) of the group $L{\bf N}%
_{+}\left( n\right) $ preserves ${{\Bbb {Y}}}_{n}$.

\begin{theorem}
\label{orbits}\cite{FRS,SemSev}

\begin{enumerate}
\item  The gauge action of $L{\bf N}_{+}\left( n\right) $ on ${{\Bbb {Y}}}%
_{n}$ is free.

\item  The set of companion matrices of the form (\ref{e110}) is a
cross-section of this action, so the quotient ${{\Bbb {Y}}}_{n}/L{\bf N}%
_{+}\left( n\right) $ can be identified with ${{\Bbb {M}}}_{n}$.
\end{enumerate}
\end{theorem}

The quotient ${{\Bbb {Y}}}_{n}/L{\bf N}_{+}\left( n\right) $ has a natural
description in the framework of Poisson reduction proposed in \cite
{FRS,SemSev}.

First of all note that the set $L{\frak {gl}}_{n}$ of all matrix first order
difference operators with potential of ''general form'' 
\begin{equation*}
{\Bbb A}=D-{\Bbb L},\quad {\Bbb L}\in L{\frak {gl}}_{n},
\end{equation*}
can be supplied with the structure of a Poisson manifold. Unlike the
differential case, the choice of this structure is not quite unique .
However, it may be fixed in a canonical way if we supply the gauge group
with the structure of a Poisson Lie group and require the Poisson bracket on 
${\cal M}_{n}$ to be {\em covariant} with respect to the gauge action; in
other words, the map 
\begin{equation*}
L\,GL_{n} \times L{\frak gl}_{n} \to L{\frak gl}_{n}:\quad {{\Bbb {L}}}\longmapsto {{\Bbb {%
L}}}^{\prime }=\,^{h}S{{\Bbb {L}}}S^{-1}.
\end{equation*}
has to be a Poisson mapping \cite{Sem}. The Poisson bracket on $L{\frak gl}_{n}$
is explicitly described in terms of the r-matrix $\hat{r}$ which fixes the
Poisson structure on gauge group (this r-matrix must satisfy the natural
invariance condition $\left( \hat{h}\otimes \hat{h}\right) \hat{r}=\hat{r}).$
We shall write down the corresponding formula a little later after we
introduce some necessary notations (see (\ref{e115} below). The reduction is
actually performed with respect to a subgroup of the full gauge group; the
natural consistency conditions are as follows:

First, the    invariants of the gauge action of $L{\bf N}_{+}\left( n\right) 
$ on $L{\frak gl}_{n}$ must form a Lie subalgebra  ${\cal I}$ with respect
to the Poisson bracket on $L{\frak gl}_{n}$; (in that case the subgroup $L%
{\bf N}_{+}\left( n\right) $ is said to be an admissible subgroup of the
full gauge group (regarded as a Poisson Lie group)). Second, the constraints
defining the submanifold ${{\Bbb {Y}}}_{n}\subset L{\frak gl}_{n}$ must
generate a Poisson ideal in ${\cal I}.$

The latter condition means that the Poisson brackets of the constraints
vanish on the constraints surface ${{\Bbb {Y}}}_{n}$, i.e., the constraints
are of the {\em first class, }according to Dirac. Both conditions impose
restrictions on the choice of the initial r-matrix (which eventually allow
to fix it completely).

To describe the relevant Poisson brackets explicitly let us fix the
following notation.

By definition, a functional $\widehat{\varphi }\in Fun\left( L{\frak gl}%
_{n}\right) $ is said to be smooth if for any ${{\Bbb {L}}}\in L{\frak gl}%
_{n}$ there exists an element $d\widehat{\varphi }\left( {{\Bbb {L}}}\right)
\in L{\frak gl}_{n}$ called its linear gradient such that

\begin{equation*}
\left\langle d\widehat{\varphi }\left( {{\Bbb {L}}}\right) ,X\right\rangle
=\left( \frac{d}{dt}\right) _{t=0}\varphi \left( {{\Bbb {L}}}+tX\right)
,\quad \forall X\in L{\frak gl}_{n}.
\end{equation*}
In applications, various functionals may be defined only on an affine
subspace of $L{\frak gl}_{n};$ in that case the choice of the gradient (when
it exists) is not unique (however, a canonical choice is frequently
possible).

For a smooth functional $\widehat{\varphi }$ we define its left and right
gradients by $\nabla \widehat{\varphi }\left( {{\Bbb {L}}}\right) ={{\Bbb {L}%
}}d\widehat{\varphi }\left( {{\Bbb {L}}}\right) $ and $\nabla ^{\prime }%
\widehat{\varphi }\left( {{\Bbb {L}}}\right) =d\widehat{\varphi }\left( {%
{\Bbb {L}}}\right) {{\Bbb {L}}},$ respectively.

Let $\hat{r}\in End\left( L{\frak gl}_{n}\right) $ be a classical r-matrix.
We assume that $\hat{r}$ is skew symmetric and satisfies the modified
classical Yang-Baxter equation and, moreover, $\hat{r}\circ \hat{h}=\hat{h}%
\circ \hat{r}.$ Put $\hat{r}_{\pm }=\hat{r}\pm \frac{1}{2}id.$ The natural
Poisson bracket on $L{\frak gl}_{n}$ which is covariant
with respect to gauge transformations is given by

\begin{equation}
\left\{ \hat{\varphi},\hat{\psi}\right\} \text{=}\left\langle \left\langle
\left( 
\begin{array}{cc}
\hat{r} & -\hat{h}\hat{r}_{+} \\ 
\hat{r}_{-}\hat{h}^{-1} & -\hat{r}
\end{array}
\right) \left( 
\begin{array}{l}
\nabla \hat{\varphi} \\ 
\nabla ^{\prime }\hat{\varphi}
\end{array}
\right) ,\left( 
\begin{array}{l}
\nabla \hat{\psi} \\ 
\nabla ^{\prime }\hat{\psi}
\end{array}
\right) \right\rangle \right\rangle  \label{e115}
\end{equation}

\noindent (see \cite{Sem}, where this formula is derived from the theory of
the so called {\em twisted Heisenberg double}). In the present context we
need to choose $\hat{r}$, so as to assure the admissibility of $L{\bf N}%
_{+}\left( n\right) $. The admissibility criterion may be found in \cite{Sem}%
; here we only describe the r-matrices which make $L{\bf N}_{+}\left(
n\right) $ admissible.

Let ${\cal P}_{+},{\cal P}_{-,}{\cal P}_{0}$ be the projection operators
onto $L{{\frak {n}}}_{+}\left( n\right) ,L{{\frak {n}}}_{-}\left( n\right) ,L%
{{\frak {h}}}_{n},$ respectively. We denote by ${\cal P}_{00}$ the
orthogonal projection operator onto the one-dimensional subspace ${{\Bbb {C}}%
}\cdot 1\subset L{{\frak {h}}}_{n}$ and put ${\cal P}_{0}^{\prime }={\cal P}%
_{0}-{\cal P}_{00}.$

Put 
\begin{equation}
\hat{r}=\frac{1}{2}\left( {\cal P}_{+}-{\cal P}_{-}\right) +\hat{r}^{0}{\cal %
P}_{0},  \label{e114}
\end{equation}
where $\hat{r}^{0}\in End\left( L{{\frak {h}}}_{n}\right) ,$ $\hat{r}%
^{0}=-\left( \hat{r}^{0}\right) ^{\ast }.$

It may be shown that with this choice of $\hat{r}$ the gauge action of $L%
{\bf N}_{+}\left( n\right) $ is admissible.

The remaining freedom in choice of $\hat{r}$ may be (almost) eliminated when
we impose our second condition. Namely, let $\hat{\tau}_{n}\in End\left( L{%
{\frak {h}}}_{n}\right) $ be the operator acting by the cyclic permutation
of matrix elements: 
\begin{equation}
\hat{\tau}_{n}{\rm diag}\left( f_{0},f_{1},\ldots ,f_{n-1}\right) ={\rm diag}%
\left( f_{1},\ldots ,f_{n-1},f_{0}\right) .  \label{e124}
\end{equation}
Let $U_{n}\subset L{{\frak {h}}}_{n}$ be the subspace of matrices of the
form 
\begin{equation}
U_{n}=\left\{ {\rm diag}\left( f_{0}\left( z\right) ,f_{0}\left(
q^{-1}z\right) ,\ldots ,f_{0}\left( q^{-\left( n-1\right) }z\right) \right)
,\quad f_{0}\in {{\Bbb {C}}}\left( \left( z^{-1}\right) \right) \right\} .
\label{e124p1}
\end{equation}

\begin{theorem}
\label{tt12}

The Poisson bracket of the form (\ref{e115}) admits the q-deformed
DS-reduction if and only if the corresponding operator $\hat{r}^{0}$ has the
form: 
\begin{equation}
\hat{r}^{0}=\frac{1}{2}\frac{1+\hat{h}\hat{\tau}_{n}}{1-\hat{h}\hat{\tau}_{n}%
}{\cal P}_{0}^{\prime }+\Delta {\cal P}_{U_{n}}+\hat{\alpha}-\hat{\alpha}%
^{*},  \label{e125}
\end{equation}
where ${\cal P}_{U_{n}}$ is the orthogonal projection operator onto $U_{n},$ 
$\Delta $ is a skew symmetric operator in $U_{n}$ and $\hat{\alpha}$ is a
linear operator of the form 
\begin{equation}
\hat{\alpha}\left( f\right) =\left\langle f,\alpha \right\rangle \cdot
1\subset L{{\frak {h}}}_{n},\quad \alpha \in U_{n}^{\bot }.  \label{e126}
\end{equation}
\end{theorem}

\begin{remark}
\label{rr11}The operator $1-\hat{h}\hat{\tau}_{n}$ is not invertible in $L{%
{\frak {h}}}_{n}.$ It is easy to see that ${\rm Ker}\left( 1-\hat{h}\hat{\tau%
}_{n}\right) ={{\Bbb {C}}}\cdot 1,$ ${\rm Im}\left( 1-\hat{h}\hat{\tau}%
_{n}\right) =L{{\frak {h}}}_{n}^{\prime },$ where 
\begin{equation}
L{{\frak {h}}}_{n}^{\prime }\equiv L{{\frak {h}}}_{n}\ominus \left\{ {{\Bbb {%
C}}}\cdot 1\right\} ={\cal P}_{0}^{\prime }\left( L{{\frak {h}}}_{n}\right) .
\label{e126p1}
\end{equation}
Hence we may define a regularized inverse operator $\left( 1-\hat{h}\hat{\tau%
}_{n}\right) ^{-1}:L{{\frak {h}}}_{n}^{\prime }\rightarrow L{{\frak {h}}}%
_{n}^{\prime }.$ In (\ref{e125}) ${\cal P}_{0}^{\prime }$ is the projection
operator onto $L{{\frak {h}}}_{n}^{\prime },$ so this expression is
well-defined.
\end{remark}

\vspace{1pt}

\begin{remark}
\label{rr12} The bracket constructed in \cite{FRS, SemSev} corresponds to 
\begin{equation}
\hat{r}_{0,n}^{0}\equiv \frac{1}{2}\frac{1+\hat{h}\hat{\tau}_{n}}{1-\hat{h}%
\hat{\tau}_{n}}{\cal P}_{0}^{\prime }.  \label{e127}
\end{equation}
\end{remark}

\vspace{1pt}

{\em Proof. }Denote by $V_{n}\subset L{{\frak {h}}}_{n}$ the space of
matrices of the form: 
\begin{equation*}
V_{n}=\left\{ {\rm diag}\left( 0,\ast ,\ldots ,\ast \right) \right\} .
\end{equation*}

\begin{lemma}
\label{pp11} Condition 2) above is equivalent to the following equation
for $\hat{r}^{0}$: 
\begin{equation}
\hat{r}^{0}\left( 1-\hat{h}\hat{\tau}_{n}\right) f=
\frac{1}{2}\left( 1+\hat{h}\hat{\tau}_{n}\right) f+\tilde{\alpha}\left(
f\right), \quad \forall f\in V_{n},  \label{e130}
\end{equation}
where $\tilde{\alpha}$ is a linear operator in $L{{\frak {h}}}_{n}$ with $%
{\rm Im}\tilde{\alpha}\subset {{\Bbb {C}}}\cdot 1\subset L{{\frak {h}}}_{n}.$
\end{lemma}

We omit the proof, since it is similar to the proof of proposition \ref{pp41}
below.

Let us define the following subspaces in $L{{\frak {h}}}_{n}$: 
\begin{equation}
\begin{array}{l}
U_{n}^{\prime }=L{{\frak {h}}}_{n}^{\prime }\cap U_{n}; \\ 
{\rm Im}V_{n}=\left( 1-\hat{h}\hat{\tau}_{n}\right) V_{n}.
\end{array}
\label{e123}
\end{equation}
It is easy to see that 
\begin{equation}
L{{\frak {h}}}_{n}=U_{n}\oplus {\rm Im}V_{n}={{\Bbb {C}}}\cdot 1\oplus
U_{n}^{\prime }\oplus {\rm Im}V_{n}.  \label{e123p1}
\end{equation}
The following lemma finishes our arguments:

\begin{lemma}
\label{lll2} Any skew-symmetric operator $\hat{r}^{0}\in EndL{{\frak {h}}}%
_{n}$ satisfying (\ref{e130}) has the form (\ref{e125}).
\end{lemma}

$\odot$ Put ${\cal D}=\hat{r}^{0}-\hat{r}_{0,n}^{0}.$ It is evident that $%
\hat{r}_{0,n}^{0}$ satisfies (\ref{e130}) with $\tilde{\alpha}=0$ (because ${%
{\Bbb {C}}}\cdot 1\perp {\rm Im}V_{n}$), hence ${\cal D}$ satisfies 
\begin{equation}
\tilde{\alpha}\left( f\right) ={\cal D}\left( 1-\hat{h}\hat{\tau}_{n}\right)
f,\quad \forall f\in V_{n}.  \label{e132}
\end{equation}
We put $\hat{\alpha}=\tilde{\alpha}\left( 1-\hat{h}\hat{\tau}_{n}\right)
^{-1}$ and rewrite (\ref{e132}) as 
\begin{equation}
\left( {\cal D}-\hat{\alpha}\right) \bar{f}=0,\quad \forall \bar{f}\in {\rm %
Im}V_{n}.  \label{e133}
\end{equation}
The operator $\hat{\alpha}$ is defined only on the subspace ${\rm Im}%
V_{n}\approx \left( {\rm Im}V_n\right)^\ast,$ hence it can be written in the
form $\hat{\alpha}= \left\langle \cdot ,\alpha \right\rangle \cdot 1$ with
some $\alpha \in {\rm Im}V_{n}=U_{n}^{\bot }.$

Skew-symmetry of ${\cal D}$ implies that it has the following block form
with respect to the orthogonal decomposition $L{{\frak {h}}}_{n}={{\Bbb {C}}}%
\cdot 1\oplus U_{n}^{\prime }\oplus {\rm Im}V_{n}$: 
\begin{equation}
\begin{array}{ll}
\begin{array}{c}
\begin{array}{rccc}
\makebox[1.5cm]{\ } & \makebox[1.5cm]{$ {{\Bbb {C}}}\cdot 1 $} & %
\makebox[2cm]{$ U_{n}^{\prime }$} & \makebox[2.5cm]{$ {\rm Im}V_{n}$}
\end{array}
\\ 
\begin{array}{r|c|c|c|}
\cline{2-4}
\makebox[1.5cm]{$\hfill {{\Bbb {C}}}\cdot 1 $} & \makebox[1.5cm]{$ 0$} & %
\makebox[2cm]{$ \beta$} & \makebox[2.5cm]{$\gamma$} \\ \cline{2-4}
U_{n}^{\prime } & -\beta ^{*} & a & b \\ \cline{2-4}
{\rm Im}V_{n} & -\gamma ^{*} & -b^{*} & d \\ \cline{2-4}
\end{array}
\end{array}
& 
\begin{array}[t]{r}
a=-a^{*}, \\ 
d=-d^{*}.
\end{array}
\end{array}
\label{e134}
\end{equation}
The equation (\ref{e133}) implies $\gamma =\hat{\alpha},$ $b=d=0.$ Put 
\begin{equation*}
\Delta =\left( 
\begin{array}{cc}
0 & \beta \\ 
-\beta ^{*} & a
\end{array}
\right) \in End\left( U_{n}\right) ,
\end{equation*}
then ${\cal D}=\Delta {\cal P}_{U_{m}}+\hat{\alpha}-\hat{\alpha}^{*},$ as
desired. $\blacksquare $

\begin{remark}
\label{rr13} It is easy to see that the term $\hat{\alpha}-\hat{\alpha}^{*}$
does not affect the value of the reduced bracket.

Indeed, let us denote 
\begin{equation}
Z_{\hat{\varphi}}=\,^{h^{-1}}\nabla \hat{\varphi}-\nabla ^{\prime }\hat{%
\varphi},\quad \bar{Z}_{\hat{\varphi}}=\,^{h^{-1}}\nabla \hat{\varphi}%
+\nabla ^{\prime }\hat{\varphi}.  \label{e116}
\end{equation}
The bracket (\ref{e115}) may be written as 
\begin{equation}
\left\{ \widehat{\varphi },\widehat{\psi }\right\} =\left\langle Z_{\widehat{%
\varphi }},\tfrac{1}{2}\bar{Z}_{\widehat{\psi }}-\hat{r}Z_{\widehat{\psi }%
}\right\rangle .  \label{e117}
\end{equation}
The formula (\ref{e117}) implies that the contribution $p_{\alpha }$ of the
term $\hat{\alpha}-\hat{\alpha}^{*}$ into the bracket is given by 
\begin{equation*}
p_{\alpha }=\left\langle \hat{\alpha}\left( Z_{\hat{\varphi}}\right) ,Z_{%
\hat{\psi}}\right\rangle -\left\langle \hat{\alpha}\left( Z_{\hat{\psi}%
}\right) ,Z_{\hat{\varphi}}\right\rangle =\hat{\alpha}\left( Z_{\hat{\varphi}%
}\right) \mathop{Tr}Z_{\hat{\psi}}-\hat{\alpha}\left( Z_{\hat{\psi}}\right) %
\mathop{Tr}Z_{\hat{\varphi}}=0,
\end{equation*}
since $\mathop{Tr}Z_{\hat{\psi}}=\mathop{Tr}Z_{\hat{\varphi}}=0,$ due to
invariance of the inner product.

Below we always put $\hat{\alpha}=0.$
\end{remark}

Now we calculate the $U_{n}$-block of the r-matrix $\hat{r}_{0,n}^{0}.$ This
auxiliary result will be used in section 6 to compute the quotient bracket
obtained via DS-reduction from the algebra of complex size matrices.

\begin{theorem}
\label{tt13} 
\begin{equation}
\left\langle \hat{r}_{0,n}^{0}\bar{f},\bar{g}\right\rangle =\left\langle 
\frac{n}{2}\frac{1+\hat{h}^{n}}{1-\hat{h}^{n}}{\cal P}_{0}^{\prime }\bar{f},%
\bar{g}\right\rangle ,\quad \forall \bar{f},\bar{g}\in U_{n}  \label{e138}
\end{equation}
\end{theorem}

{\em Proof. }To calculate the l.h.s. we expand $\bar{f},\bar{g}$ with
respect to the basis of eigenfunctions of the operator $\hat{h}\hat{\tau}%
_{n}.$

\begin{lemma}
\label{ll12}The eigenfunctions of the operator $\hat{h}\hat{\tau}_{n}$ are 
\begin{equation}
E_{m,\alpha }=z^{m}{\bf e}_{\alpha },\quad m\in {{\Bbb {Z}}}\;,\quad \alpha
=0,\ldots ,n-1,  \label{e138p1}
\end{equation}
where 
\begin{equation}
{\bf e}_{\alpha }={\bf diag}\left( 1,\omega ^{\alpha },\ldots ,\omega
^{\left( n-1\right) \alpha }\right) ,\quad \omega =e^{\frac{2\pi i}{n}}.
\label{e138p2}
\end{equation}
The corresponding eigenvalues $\xi _{m,\alpha }$ are equal to 
\begin{equation}
\xi _{m,\alpha }=q^{m}\omega ^{\alpha }.  \label{e138p3}
\end{equation}
The eigenfunctions satisfy 
\begin{equation}
\left\langle E_{m,\alpha },E_{l,\beta }\right\rangle =n\delta _{m,-l}\cdot
\left\{ 
\begin{array}{l}
1,\quad \alpha =-\beta \mathop{mod}n, \\ 
0,\quad \text{in other cases},
\end{array}
\right.   \label{e138p4}
\end{equation}
and form a basis in $L{{\frak {h}}}_{n}.$
\end{lemma}

We shall denote by $\sum_{m,\alpha }^{\prime }$ the sum over all pairs $%
\left( m,\alpha \right) \neq \left( 0,0\right) ,$ $m\in {{\Bbb {Z}}},$ $%
\alpha =0,\ldots ,n-1.$ Note that in the expansion of $\bar{f}$ with respect
to the eigenbasis $E_{m,\alpha }$ the $E_{0,0}$-component is annihilated by $%
{\cal P}_{0}^{\prime }$, hence 
\begin{equation}
\left\langle \hat{r}_{0,n}^{0}\bar{f},\bar{g}\right\rangle = \mathop{{\sum}'}%
_{m,\alpha }\frac{1}{n}\frac{1+q^{m}\omega ^{\alpha }}{1-q^{m}\omega
^{\alpha }}\left\langle \bar{f}\,,E_{-m,n-\alpha }\right\rangle \left\langle 
\bar{g},E_{m,\alpha }\right\rangle .  \label{e139}
\end{equation}
Any element of $U_{n}$ has the form $\bar{f}={\rm diag}\left( f\left(
z\right) ,f\left( q^{-1}z\right) ,\ldots ,f\left( q^{-\left( n-1\right)
}z\right) \right) ;$ we denote by $f_{m}$ the coefficient of the formal
Laurent expansion of $f\left( z\right) $ corresponding to $z^{m}$: $f\left(
z\right) =\sum\limits_{m=-\infty }^{N\left( f\right) }f_{m}z^{m}.$ It is
easy to see that 
\begin{equation}
\left\langle \bar{f}\,,E_{-m,n-\alpha }\right\rangle
=f_{m}\sum_{i=0}^{n-1}\left( q^{m}\omega ^{\alpha }\right) ^{-i}=f_{m}\frac{%
1-q^{-mn}\omega ^{-\alpha n}}{1-q^{-m}\omega ^{-\alpha }}=f_{m}\frac{%
1-q^{-mn}}{1-q^{-m}\omega ^{-\alpha }}.  \label{e141}
\end{equation}
Substituting (\ref{e141}) and a similar expression for $\left\langle \bar{g}%
,E_{m,\alpha }\right\rangle $ into (\ref{e139}) we find: 
\begin{equation}
\left\langle \hat{r}_{0,n}^{0}\bar{f},\bar{g}\right\rangle =\frac{1}{2}%
\mathop{{\sum}'}_{m,\alpha }f_{m}g_{-m}\left( 1-q^{mn}\right) ^{2}q^{-mn}%
\frac{1}{n}\frac{1+q^{m}\omega ^{\alpha }}{\left( 1-q^{m}\omega ^{\alpha
}\right) ^{3}}q^{m}\omega ^{\alpha }.  \label{e142}
\end{equation}
In the sum above the terms corresponding to $m=0$ vanish due to the
multiplier $1-q^{mn}.$

\begin{lemma}
\label{ll13} 
\begin{equation}
S_{0}\equiv \sum_{\alpha =0}^{n-1}\frac{1}{n}\frac{1+q^{m}\omega ^{\alpha }}{%
\left( 1-q^{m}\omega ^{\alpha }\right) ^{3}}q^{m}\omega ^{\alpha
}=n^{2}q^{mn}\frac{1+q^{mn}}{\left( 1-q^{mn}\right) ^{3}},\quad m\neq 0.
\label{e143}
\end{equation}
\end{lemma}

$\odot $ Note that 
\begin{equation*}
z\frac{1+z}{\left( 1-z\right) ^{3}}=\sum_{k=1}^{\infty }k^{2}z^{k},\quad
|z|<1,
\end{equation*}
hence 
\begin{equation*}
S_{0}=\sum_{k=1}^{\infty }q^{mk}k^{2}\frac{1}{n}\sum_{\alpha =0}^{n-1}\omega
^{\alpha k}.
\end{equation*}
But 
\begin{equation*}
\frac{1}{n}\sum_{\alpha =0}^{n-1}\omega ^{\alpha k}=\left\{ 
\begin{array}{l}
0,\quad k\neq jn,\quad j\in {{\Bbb {Z}}}, \\ 
1,\quad k=jn,
\end{array}
\right.
\end{equation*}
therefore 
\begin{equation*}
S_{0}=\sum_{j=1}^{\infty }q^{mnj}\left( nj\right) ^{2}=n^{2}q^{mn}\frac{%
1+q^{mn}}{\left( 1-q^{mn}\right) ^{3}}.\qquad \blacksquare
\end{equation*}
Using the lemma we find from (\ref{e142}): 
\begin{eqnarray*}
\left\langle \hat{r}_{0,n}^{0}\bar{f},\bar{g}\right\rangle &=&\frac{n^{2}}{2}%
\sum_{m\neq 0}f_{m}g_{-m}\frac{1+q^{mn}}{1-q^{mn}}=\frac{n^{2}}{2}%
\left\langle \left[ \frac{1+\hat{h}^{n}}{1-\hat{h}^{n}}\left( 1-\mathop{Res}%
\right) \right] f,g\right\rangle _{{{\Bbb {C}}}} \\
&=&\left\langle \frac{n}{2}\frac{1+\hat{h}^{n}}{1-\hat{h}^{n}}{\cal P}%
_{0}^{\prime }\bar{f},\bar{g}\right\rangle _{L{{\frak {h}}}_{n}},
\end{eqnarray*}
as desired. $\blacksquare $

\subsection{Explicit formula for the quotient bracket.}

As mentioned above, the quotient ${{\Bbb {Y}}}_{n}/LN_{+}\left( n\right) $
can be identified with the space ${{\Bbb {M}}}_{n}$ of scalar q-difference
operators of $n$-th order. To describe the quotient bracket we shall
consider ${{\Bbb {M}}}_{n}$ as an affine subspace in the algebra $\Psi {\bf D%
}_{q}$ of q-pseudodifference symbols. By definition, $\Psi {\bf D}_{q}$
consists of formal series of the form 
\begin{equation}
A=\sum_{i=-\infty }^{N\left( A\right) }a_{i}\left( z\right) D^{i},\quad
a_{i}\in {{\Bbb {C}}}\left( \left( z^{-1}\right) \right)  \label{e150}
\end{equation}
with the commutation relation 
\begin{equation}
D\cdot a=\,^{h}a\cdot D.  \label{e151}
\end{equation}

As a linear space, $\Psi {\bf D}_{q}$ is a direct sum of three subalgebras, 
\begin{align}
J_{+}& =\left\{ A\in \Psi {\bf D}_{q}|\;A=\sum_{i=1}^{N\left( A\right)
}a_{i}\left( z\right) D^{i},\quad a_{i}\in {{\Bbb {C}}}\left( \left(
z^{-1}\right) \right) \right\} ,  \label{e152} \\
J_{0}& ={{\Bbb {C}}}\left( \left( z^{-1}\right) \right) \subset \Psi {\bf D}%
_{q},  \label{e153} \\
J_{-}& =\left\{ A\in \Psi {\bf D}_{q}|\;A=\sum_{i=1}^{\infty }a_{i}\left(
z\right) D^{-i},\quad a_{i}\in {{\Bbb {C}}}\left( \left( z^{-1}\right)
\right) \right\} .  \label{e154}
\end{align}
Clearly, $J_{0}$ normalizes $J_{\pm }$ and hence $J_{\left( \pm \right)
}=J_{\pm }+J_{0}$ is also a subalgebra. Let $P_{\pm },P_{0}$ be the
associated projection operators which project $\Psi {\bf D}_{q}$ onto $%
J_{\pm },J_{0},$ respectively, parallel to the complement. Put $P_{\left(
\pm \right) }=P_{\pm }+P_{0}.$ For $A\in \Psi {\bf D}_{q}$ set $A_{\pm }=
P_{\pm }A,\quad A_{\left( \pm \right) }=P_{\left( \pm \right) }A.$

We define the residue of a q-pseudodifference operator $A$ by 
\begin{equation*}
{\rm {\mathop{Res}}}A=A_{0}=P_{0}A.
\end{equation*}
It is easy to see that the formal trace defined by 
\begin{equation}
Tr\;A=\int {\rm {\mathop{Res}}}Adz/z  \label{e155}
\end{equation}
satisfies the natural condition 
\begin{equation*}
Tr\;AB=Tr\;BA,\quad A,B\in \Psi {\bf D}_{q}.
\end{equation*}
We introduce an inner product in $\Psi {\bf D}_{q}$ by 
\begin{equation}
\left\langle A,B\right\rangle =Tr\;AB,\quad A,B\in \Psi {\bf D}_{q}.
\label{e156p}
\end{equation}
Clearly, this inner product is invariant and non-degenerate and the
subalgebras $J_{\pm }$ are isotropic; moreover, it sets $J_{+}$ and $J_{-}$
into duality, while $J_{0}\simeq J_{0}^{*}.$

We shall now define a class of Poisson brackets on $\Psi {\bf D}_{q}.$ The
natural algebra of observables $Fun\left( \Psi {\bf D}_{q}\right) $ in the
present case is generated by 'elementary' functionals which assign to a
pseudodifference operator $A$ the formal integrals of its coefficients, 
\begin{equation*}
\zeta _{i}^{j}\left( A\right) =\mathop{Tr}\left( z^{-j}AD^{-i}\right) .
\end{equation*}
As compared to the case of differential operators, the definition of a
quadratic Poisson bracket in the difference case is not quite
straightforward; the point is that the 'naive' bracket defined by analogy
with the differential case is not compatible with the natural normalization
condition for difference operators (highest coefficient is set to one); an
easy scrutiny shows that the source of the trouble lies in the $J_{0}$%
-component in the expansion 
\begin{equation*}
\Psi {\bf D}_{q}=J_{+}\dot{+}J_{0}\dot{+}J_{-}.
\end{equation*}
To avoid this difficulty we are bound to consider a more general class of
quadratic Poisson bracket which mix together both left and right gradients
of functions\footnote{%
This class of Poisson brackets naturally arises in the theory of Poisson Lie
groups, cf. \cite{Monod}, \cite{FM}, \cite{LP}.}. (Recall that
in the Gelfand-Dickey case left and right gradients are coupled only to the
gradients of the same chirality.)

For a smooth functional $\varphi $ let us write 
\begin{equation*}
D\varphi =\left( 
\begin{array}{c}
\nabla \varphi \\ 
\nabla ^{\prime }\varphi
\end{array}
\right) .
\end{equation*}

Let us consider quadratic Poisson brackets on $\Psi {\bf D}_{q}$ of the
following form:

\begin{equation}
\left\{ \varphi ,\psi \right\} =\left\langle \left\langle \left( 
\begin{array}{cc}
R+aP_{0} & bP_{0} \\ 
cP_{0} & R+dP_{0}
\end{array}
\right) D\varphi ,D\psi \right\rangle \right\rangle ,  \label{e156}
\end{equation}
where $R=\frac{1}{2}\left( P_{+}-P_{-}\right) $ and $a,b,c,d$ are linear
operators acting in $J_{0}$ satisfying 
\begin{equation*}
a=-a^{*},\quad d=-d^{*},\quad c^{*}=b.
\end{equation*}
In other words, the bracket (\ref{e156}) differs from the naive
Gelfand-Dickey bracket by a 'perturbation term' which is acting only on the $%
J_{0}$-components of the gradients. This bracket satisfies the Jacobi
identity for any choice of $a,b,c,d$. Note that different $a,b,c,d$ may give
rise to the same bracket. More precisely, we have the following

\begin{lemma}
\label{ll14} Let $f,g,h,k$ be linear operators in $J_{0}$ with images in the
subspace of constants ${{\Bbb {C}}}\cdot 1\subset J_{0}.$ The r-matrices $%
{\cal R}$ and ${\cal R}^{\prime }={\cal R}+\Theta $ where 
\begin{equation*}
\Theta =\left( 
\begin{array}{cc}
h-k^{*} & f+k^{*} \\ 
h+g^{*} & -g^{*}+f
\end{array}
\right) ,
\end{equation*}
define the same Poisson bracket.
\end{lemma}

Up to this ambiguity, the unique choice of the coefficients $a,b,c,d$ is
assured by the condition that the set ${\Bbb M}_{n}$ of difference operators
with normalized highest coefficient is a Poisson submanifold with respect to
the Poisson structure (\ref{e156}) and that, moreover, formal spectral
invariants of difference operators give rise to Lax equations of standard
commutator form. More precisely, we have the following theorem (see \cite
{PirSem}):

\begin{theorem}
\label{tt16}There exists a unique Poisson bracket of the form (\ref{e156})
on $\Psi {\bf D}_{q}$ such that

\begin{itemize}
\item[1)]  the affine subspace ${{\Bbb {M}}}_{n}$ is a Poisson submanifold;

\item[2)]  Formal spectral invariants $H_{m}=\frac{n}{m}\mathop{Tr}L^{\frac{m%
}{n}},\quad m\in {{\Bbb {N}}},$are in involution.
\end{itemize}

This bracket is given by \footnote{%
For difference Lax equations on the lattice Poisson bracket \ref{e} was also
introduced in \cite{Oew}.}
\begin{equation}
\left\{ \varphi ,\psi \right\} _{n}^{0}=\left\langle \left\langle \left( 
\begin{array}{cc}
R+\left( \frac{1}{2}\frac{1+\hat{h}^{n}}{1-\hat{h}^{n}}\right) P_{0}^{\prime
} & -\left( \frac{\hat{h}^{n}}{1-\hat{h}^{n}}\right) P_{0}^{\prime } \\ 
\left( \frac{1}{1-\hat{h}^{n}}\right) P_{0}^{\prime } & R-\left( \frac{1}{2}%
\frac{1+\hat{h}^{n}}{1-\hat{h}^{n}}\right) P_{0}^{\prime }
\end{array}
\right) D\varphi ,D\psi \right\rangle \right\rangle .  \label{e}
\end{equation}
\end{theorem}

\begin{remark}
It is easy to see that  the  involutivity condition (2)  is equivalent to
the following simple constraint: 
\begin{equation}
a+b=c+d.  \label{e158}
\end{equation}
\end{remark}

The quotient bracket on the set of q-difference operators which is obtained
via the q-DS reduction differs from the above formula by an additional term
which reflects the remaining freedom in the choice of the classical r-matrix
on $L{\frak gl}_{n}$ which is compatible with the reduction; namely:

\begin{theorem}
\label{tt14} Let 
\begin{equation}
\hat{r}_{\Delta ,n}^{0}=\frac{1}{2}\frac{1+\hat{h}\hat{\tau}_{n}}{1-\hat{h}%
\hat{\tau}_{n}}{\cal P}_{0}^{\prime }+n\Delta {\cal P}_{U_{n}}.  \label{e159}
\end{equation}
where ${\cal P}_{U_{n}}$ is the orthogonal projection operator onto $U_{n}$
and $\Delta $ is a skew symmetric operator in $U_{n}$ commuting with $\hat{h}%
.$ Let $\hat{r}_{\Delta ,n}=\frac{1}{2}\left( {\cal P}_{+}-{\cal P}%
_{-}\right) +\hat{r}_{\Delta ,n}^{0}{\cal P}_{0}$.

The Poisson bracket $\left\{ \cdot ,\cdot \right\} _{\Delta ,n}$ on $L{%
{\frak {gl}}}_{n}$ defined by 
\begin{equation}
\left\{ \hat{\varphi},\hat{\psi}\right\} _{\Delta ,n}\text{=}\left\langle
\left\langle \left( 
\begin{array}{cc}
\hat{r}_{\Delta ,n} & -\hat{h}\left( \hat{r}_{\Delta ,n}\right) _{+} \\ 
\left( \hat{r}_{\Delta ,n}\right) _{-}\hat{h}^{-1} & -\hat{r}_{\Delta ,n}
\end{array}
\right) \left( 
\begin{array}{l}
\nabla \hat{\varphi} \\ 
\nabla ^{\prime }\hat{\varphi}
\end{array}
\right) ,\left( 
\begin{array}{l}
\nabla \hat{\psi} \\ 
\nabla ^{\prime }\hat{\psi}
\end{array}
\right) \right\rangle \right\rangle   \label{e160}
\end{equation}
gives rise via DS-reduction to the following bracket on ${{\Bbb {M}}}_{n}$: 
\begin{equation}
\begin{array}{l}
\left\{ \varphi ,\psi \right\} _{n}^{\Delta }=\left\langle \left\langle
\left( 
\begin{array}{cc}
R+\left( \frac{1}{2}\frac{1+\hat{h}^{n}}{1-\hat{h}^{n}}+\Delta \right)
P_{0}^{\prime } & -\left( \frac{1}{1-\hat{h}^{n}}+\Delta \right) \hat{h}%
^{n}P_{0}^{\prime } \\ 
\left( \frac{\hat{h}^{n}}{1-\hat{h}^{n}}+\Delta \right) \hat{h}%
^{-n}P_{0}^{\prime } & R-\left( \frac{1}{2}\frac{1+\hat{h}^{n}}{1-\hat{h}^{n}%
}+\Delta \right) P_{0}^{\prime }
\end{array}
\right) D\varphi ,D\psi \right\rangle \right\rangle ; \\ 
\quad  \\ 
\quad 
\end{array}
\label{e161}
\end{equation}
(here we have identified $EndU_{n}$ and $End{{\Bbb {C}}((}z^{-1})).$ )
\end{theorem}

The remaining ambiguity in the choice of r-matrix may be removed if  we
impose the involutivity condition.

\begin{theorem}
\label{tt15} The only one of brackets $\left\{ \cdot ,\cdot \right\}
_{\Delta ,n}$ which gives rise to a Poisson bracket on ${{\Bbb {M}}}_{n}$
satisfying the involutivity condition (\ref{e158}) is $\left\{ \cdot ,\cdot
\right\} _{0,n}.$
\end{theorem}

{\em Proof.} In the class of the brackets (\ref{e161}) only the bracket $%
\left\{ \cdot ,\cdot \right\} _{n}^{0}$ satisfies this condition. Indeed, (%
\ref{e158}) implies that $\Delta \left( 2-\hat{h}^{n}-\hat{h}^{-n}\right) =0$
but the operator $2-\hat{h}^{n}-\hat{h}^{-n}$ is invertible.$\blacksquare $

{\em Proof of theorem \ref{tt14}. }For $\Delta =0$ this theorem has been
proved in \cite{PirSem}$,$ hence we need to calculate only the contribution
of the term $n\Delta {\cal P}_{U_{n}}.$ Let $\varphi $ be a smooth
functional on ${{\Bbb {M}}}_{n},$ $\hat{\varphi}$ the corresponding $%
LN_{+}\left( n\right) $-invariant functional on $L{\frak gl}_{n}.$ To fix
their gradients we shall assume that 
\begin{equation}
d\varphi =\sum_{i=0}^{n-1}f_{i}D^{-i},\quad f_{i}\in {{\Bbb {C}}}((z^{-1})),
\label{e163}
\end{equation}
and $d\hat{\varphi}\in L{{\frak {b}}}_{-}\left( n\right) .$ For $L\in {{\Bbb 
{M}}}_{n}$ let us denote by ${{\Bbb {L}}}\in L{\frak gl}_{n}$ the
corresponding companion matrix.

Let us denote by $J_{\Delta }\left( \hat{\varphi},\hat{\psi}\right) $ the
contribution of $n\Delta {\cal P}_{U_{n}}$ to the bracket (\ref{e160}). From
(\ref{e117}) it follows that 
\begin{equation}
J_{\Delta }\left( \hat{\varphi},\hat{\psi}\right) =\left\langle n\Delta 
{\cal P}_{U_{n}}Z_{\hat{\varphi}}^{0},Z_{\hat{\psi}}^{0}\right\rangle ,
\label{e164}
\end{equation}
where $Z_{\hat{\varphi}}^{0}\equiv {\cal P}_{0}Z_{\hat{\varphi}},\quad Z_{%
\hat{\psi}}^{0}\equiv {\cal P}_{0}Z_{\hat{\psi}}.$

\begin{lemma}
\label{ll17}We have 
\begin{equation}
Z_{\hat{\varphi}}^{0}\left( {{\Bbb {L}}}\right) ={\rm diag}\left(
\,^{h^{-1}}P_{0}\nabla \varphi \left( L\right) ,0,\ldots ,0,-P_{0}\nabla
^{\prime }\varphi \left( L\right) \right) .  \label{e165}
\end{equation}
\end{lemma}

\begin{lemma}
\label{ll18}The projection operator ${\cal P}_{U_{n}}$ is given by 
\begin{equation}
{\cal P}_{U_{n}}\cdot {\rm diag}\left( F_{0}\left( z\right) ,\ldots
F_{n-1}\left( z\right) \right) ={\rm diag}\left( f_{0}\left( z\right)
,f_{0}\left( q^{-1}z\right) ,\ldots ,f_{0}\left( q^{-\left( n-1\right)
}z\right) \right) ,  \label{e166}
\end{equation}
where 
\begin{equation}
f_{0}\left( z\right) =\frac{1}{n}\sum_{i=0}^{n-1}F_{i}\left( q^{i}z\right) .
\label{e167}
\end{equation}
\end{lemma}

Using these lemmas we find 
\begin{equation*}
{\cal P}_{U_{n}}Z_{\hat{\varphi}}^{0}=\frac{1}{n}\left(
\,^{h^{-1}}P_{0}\nabla \varphi -\,^{h^{n-1}}P_{0}\nabla ^{\prime }\varphi
\right) .
\end{equation*}

Substituting this into (\ref{e164}) and taking into account the invariance
of the inner product we obtain 
\begin{equation}
J_{\Delta }\left( \hat{\varphi},\hat{\psi}\right) =\left\langle \left\langle
\left( 
\begin{array}{cc}
\Delta P_{0} & -\Delta \hat{h}^{n}P_{0} \\ 
\Delta \hat{h}^{-n}P_{0} & -\Delta P_{0}
\end{array}
\right) D\varphi ,D\psi \right\rangle \right\rangle .  \label{e168}
\end{equation}
But $\Delta $ is skew-symmetric, hence it annihilates the one-dimensional
subspace ${{\Bbb {C}}}\cdot 1\subset {{\Bbb {C}}}\left( \left( z^{-1}\right)
\right) $ and we may replace $P_{0}$ by $P_{0}^{\prime }$ in (\ref{e168}).
$\blacksquare $

\section{Algebras ${{\frak {gl}}}_{q}^{\lambda }$ of complex size matrices
and their loop algebras.}

In this section we construct an algebra ${{\frak {gl}}}_{q}$ consisting of ${%
{\frak {gl}}}_{\infty }$-matrices whose matrix elements are holomorphic
functions of special form. Then we define a trace functional on ${{\frak {gl}%
}}_{q}$ with values in the space $Hol\left( {{\Bbb {C}}}\right) $ of
holomorphic functions; it satisfies the natural condition $\mathop{Tr}AB=%
\mathop{Tr}BA.$

For any fixed $\lambda \in {\Bbb {C}}$ the algebra ${{\frak {gl}}}%
_{q}^{\lambda }\subset {{\frak {gl}}}_{\infty }$ is the image of ${{\frak {gl%
}}}_{q}$ under the evaluation map $A\mapsto A\left( \lambda \right) .$ The
functional $\mathop{Tr}$ on ${{\frak {gl}}}_{q}$ induces a ${{\Bbb {C}}}$%
-valued trace on ${{\frak {gl}}}_{q}^{\lambda }.$ This construction is a
q-difference counterpart of the one described in \cite{KhM}; in particular, $%
{{\frak {gl}}}_{q}$ and ${{\frak {gl}}}_{q}^{\lambda }$ are some extensions
of the algebras ${{\frak {gl}}},\overline{{{\frak {gl}}}}_{\lambda }$
considered there. At the end of this section we shall describe the loop
algebras $L{{\frak {gl}}}_{q},$ $L{{\frak {gl}}}_{q}^{\lambda }.$

\subsection{${\cal A}_{0}$-functions.}

We shall describe a class ${\cal A}_{0}\subset Hol\left( {{\Bbb {C}}}\right) 
$ of holomorphic functions we shall deal with throughout this article.

By definition, ${\cal A}_{0}$ is the algebra of functions of complex
variable $w$ generated by $w,q^{w},q^{-w},$ where $q^{w}\equiv \exp \left(
w\ln q\right) .$ The elements of ${\cal A}_{0}$ will be called ${\cal A}_{0}$%
-functions. Evidently, the set of elements $\zeta _{m,n}=w^{m}q^{nw},$ $m\in 
{\Bbb Z}_{+},$ $n\in {\Bbb Z},$ is a linear basis of ${\cal A}_{0},$ i.e.
any ${\cal A}_{0}$-function $f$ can be decomposed into a {\em finite sum }%
with respect to this basis:

\begin{equation}
f\left( w\right) =\sum_{m+|n|\leq N(f)}f_{m,n}\zeta _{m,n},\quad f_{m,n}\in {%
{\Bbb {C}}}.  \label{e26}
\end{equation}
The minimal possible value of $N(f)$ in the sum (\ref{e26}) is called the
degree of $f$ and will be denoted by $\deg f.$ Note also that the set of
subspaces ${{\Bbb {C}}}\zeta _{m,n}$ defines a ${\Bbb Z}_{+}\times {\Bbb Z}$%
-grading on ${\cal A}_{0}.$

${\cal A}_{0}$-functions satisfy two important properties which will be
widely used below. The first one called {\em interpolation property} allows
to reconstruct an ${\cal A}_{0}$-functions from its values at sufficiently
large integer points:

\begin{proposition}
\label{p11}

Let $f\in {\cal A}_{0}$ and $f\left( n\right) =0$ for all sufficiently large
integer $n,$ then $f\left( w\right) \equiv 0.$
\end{proposition}

Hence if some relation for ${\cal A}_{0}$-functions holds for sufficiently
large integer values of $w,$ it holds also for all $w\in {{\Bbb {C}}}$.

The second property is given by

\begin{proposition}
\label{p11p}

For any ${\cal A}_{0}$-function $f$ and any $l\in {{\Bbb {Z}}}$ there exists
a unique ${\cal A}_{0}$-function $\tilde{F}$ which interpolates the sum $%
F\left( n\right) =\sum\limits_{i=0}^{n-1}f\left( i\right) q^{il},$ i.e., $%
\tilde{F}\left( n\right) =F\left( n\right) ,$ $n\in {{\Bbb {N}}}.$
\end{proposition}

\vspace{1pt}It will play the key role in the definition of trace as well as
in the proof of the cross-section theorem, see below.

We say that a function $f(w,t)$ is an ${\cal A}_{0}$-function of two complex
variables $w,t$ if it can be written as a finite sum 
\begin{equation*}
f(w,t)=\sum_{i=1}^{N(f)}f_{i}^{(1)}(w)\cdot f_{i}^{(2)}\left( t\right)
,\quad f_{i}^{(1)},f_{i}^{(2)}\in {\cal A}_{0}.
\end{equation*}
In other words, the space of ${\cal A}_{0}$-function of two variables is the
algebraic tensor product ${\cal A}_{0}\otimes {\cal A}_{0}.$

\begin{proposition}
\label{pp21}Let $f(w,t)$ be an ${\cal A}_{0}$-function of variable $w$ for
any fixed $t$ and a ${\cal A}_{0}$-function of variable $t$ for any fixed $%
w, $ then it is an ${\cal A}_{0}$-function of two variables.
\end{proposition}

\subsection{Algebras ${{\frak {gl}}}_{q},$ ${{\frak {gl}}}_{q}^{\lambda
} $ and trace functional.}

Let ${\frak a}$ be an associative algebra. We define ${{\frak {gl}}}_{\infty
}\left( {\frak a}\right) $ as the algebra of semi-infinite matrices $%
A=\left\{ A_{i,j}\in {\frak a}\right\} _{i,j=0,1,...}$, such that $A_{i,j}=0$
if $i-j>N\left( A\right) .$ For ${{\frak {gl}}}_{\infty }\left( {{\Bbb {C}}}%
\right) $ we write simply ${{\frak {gl}}}_{\infty }.$ Note that if ${\frak a}
$ is infinite dimensional, the algebra ${{\frak {gl}}}_{\infty }\left( 
{\frak a}\right) $ is wider than the algebraic tensor product ${{\frak {gl}}}%
_{\infty }\otimes {\frak a}.$

\begin{definition}
\label{d11} The algebra ${{\frak {gl}}}_{q}\subset {{\frak {gl}}}_{\infty
}\left( {\cal A}_{0}\right) $ consists of ${{\frak {gl}}}_{\infty }$%
-matrices $A\left( t\right) =\left\{ A_{i,j}\left( t\right) \right\} $ with
coefficients in ${\cal A}_{0}$ satisfying the following conditions:

There exists an integer $N\left( A\right) $ such that

\begin{enumerate}
\item  $A_{i,i+n}\left( t\right) =0$ if $n<-N\left( A\right) ;$

\item  for any fixed $n\geq -N\left( A\right) $ and any $i>N\left( A\right) $
$\ A_{i,i+n}\left( t\right) $ considered as a function of variables $i,t$
can be interpolated by an ${\cal A}_{0}$-function of two variables;

\item  For all integer $m>N\left( A\right) $ and $N\left( A\right) <i<m,$ $%
j\geq m$ we have $A_{i,j}\left( m\right) =0.$
\end{enumerate}

\vspace{1pt}The minimal possible value of $N(A)$ is called {\em the
regularity degree of} $A$ and will be denoted by ${\rm reg}A$

\end{definition}

In other words, condition 2 means that

i) for any fixed $t$ $A_{i,i+n}\left( t\right) $ considered as a function of 
$i$ can be interpolated by an ${\cal A}_{0}$-function;

ii) the degree of $A_{i,i+n}\left( t\right) $ considered as an ${\cal {A}_0}$%
-function of $t$ is uniformly bounded for all $i.$

Condition 3 means that  the matrix $A\left( m\right) \subset {{\frak {gl}}}%
_{\infty }$ has the form:

\begin{equation*}
A\left( m\right) =\left( 
\begin{array}{c|c}
a\in {{\frak {gl}}}_{m} & b \\ \hline
\begin{array}{c}
\quad  \\ 
\ast 
\end{array}
& 
\begin{array}{c}
\quad  \\ 
\ast 
\end{array}
\\ 
& 
\end{array}
\right) ,
\end{equation*}
where the number of non-zero rows in the right upper block $b$ does not
exceed ${\rm reg}A$ and hence is uniformly bounded  for all $m.$

We define the following ${{\Bbb {Z}}}$-grading on ${{\frak {gl}}}_{q}$: the
set of elements of level $n$ consists of matrices with only $n$-th non-zero
diagonal, i.e. $A_{i,i+k}\left( t\right) =0$ if $k\neq n.$ We will denote
by $A^{\left( n\right) }$ the ${\cal S}_{n}$-component of a matrix $A\in {%
{\frak {gl}}}_{q}.$

For any matrix $A\in {\cal S}_{n}$ we can assign a ${\cal A}_{0}$-function
of two variables. By definition of ${{\frak {gl}}}_{q},$ there exists an $%
{\cal A}_{0}$-function $f\left( w,t\right) $ which interpolates $%
A_{i,i+n}\left( t\right) $ {\em for all sufficiently large }$i$: 
\begin{equation*}
A_{i,i+n}\left( t\right) =f\left( i,t\right) ,\quad \forall i>N\left(
A\right) ,\quad \forall t\in {{\Bbb {C}}}.
\end{equation*}
{\em We will denote it by} $A\left( w,t\right) .$

So, for $A\in {{\frak {gl}}}_{q}$ its $n$-th diagonal $A^{\left( n\right)
}\in {\cal S}_{n}$ and $A^{\left( n\right) }\left( w,t\right) $ is the
corresponding ${\cal A}_{0}$-function. For positive integer $n$ these
functions satisfy the following important property which ensures the
invariance of the trace functional on ${{\frak {gl}}}_{q}$:

\begin{proposition}
\label{p12}

For any $A\in {{\frak {gl}}}_{q}$, $n\in {{\Bbb {N}}},$ we have: 
\begin{equation}
A^{\left( n\right) }\left( w,w+l\right) =0,\quad \forall l=1,\ldots ,n,\quad
\forall w\in {{\Bbb {C}}}.  \label{e21}
\end{equation}
\end{proposition}

\vspace{1pt}{\em Proof. }By definition of ${{\frak {gl}}}_{q},$ $%
A_{i,i+n}\left( m\right) =0$ if $m>N\left( A\right) ,$ $N\left( A\right)
<i<m,$ and $i+n\geq m,$ or, equivalently, if $i=m-l,$ $\forall l=1,\ldots
,n. $ But for any $i>N\left( A\right) $ $A_{i,i+n}\left( m\right) $
coincides with its interpolating ${\cal A}_{0}$-function. Hence $A^{\left(
n\right) }\left( m-l,m\right) =0$ for any integer $m>N\left( A\right) +l,$
and therefore by proposition \ref{p11} $A^{\left( n\right) }\left(
w-l,w\right) =0$ for $\forall w\in {{\Bbb {C}}},$ which is equivalent to (%
\ref{e21}). $\blacksquare $

Now we shall define the trace functional on ${{\frak {gl}}}_{q}.$ This
construction is parallel to the one described by Khesin and Malikov in \cite
{KhM} and goes back to J.Bernstein.

Let us consider the sum 
\begin{equation}
F_{A}\left( n,t\right) =\sum_{i=0}^{n-1}A_{i,i}\left( t\right) .  \label{e22}
\end{equation}
By proposition \ref{p11p}, there exists a unique ${\cal A}_{0}$-function of
two variables which coincides with $F_{A}\left( n,t\right) $ for any
sufficiently large integer $n.$ We will denote it by ${{\Bbb {D}}}_{A}\left(
w,t\right) .$ By definition, 
\begin{equation}
\left( \mathop{Tr}A\right) \left( t\right) ={{\Bbb {D}}}_{A}\left(
t,t\right) \in {\cal A}_{0}.  \label{e23}
\end{equation}

\begin{proposition}
\label{p13} 
\begin{equation*}
\mathop{Tr}AB=\mathop{Tr}BA.
\end{equation*}
\end{proposition}

\vspace{1pt}

{\em Proof. }It is evident that $\mathop{Tr}$ is consistent with the
grading, i.e. for any elements $A\in {\cal S}_{i},$ $B\in {\cal S}_{j}$ \
the trace of their product vanishes unless $i+j=0.$ Hence it is sufficient
to consider the case of $A\in {\cal S}_{k},$ $B\in {\cal S}_{-k},$ $k\in {%
{\Bbb {N}}}.$ For all sufficiently large $n\in {{\Bbb {N}}}$ we have 
\begin{eqnarray*}
{{\Bbb {D}}}_{AB}\left( n,t\right) -{{\Bbb {D}}}_{BA}\left( n,t\right)
&=&F_{AB}\left( n,t\right) -F_{AB}\left( n,t\right)
=\sum_{i=n-k}^{n-1}A(i,t)B(i+k,t) \\
&=&\sum_{j=1}^{k}A(n-j,t)B(n-j+k,t).
\end{eqnarray*}
All terms of this expression are ${\cal A}_{0}$-functions of variable $n,$
therefore, by the interpolation property, it holds for all complex values of 
$n;$ in particular, for $n=t$ we obtain 
\begin{equation*}
\mathop{Tr}AB-\mathop{Tr}BA\equiv {{\Bbb {D}}}_{AB}\left( t,t\right) -{{\Bbb 
{D}}}_{BA}\left( t,t\right) =\sum_{j=1}^{k}A(t-j,t)B(t-j+k,t).
\end{equation*}
But $A(t-j,t)=0$ for $j=1,\ldots ,k,$ by proposition \ref{p12}. $%
\blacksquare $

We use the following notation: ${{\frak {b}}}_{+}({{\frak {b}}}_{-}) \subset 
{{\frak {gl}}}_{q}$ is the subalgebra of upper (lower) triangular matrices, $%
{{\frak {n}}}_{+}({{\frak {n}}}_{-})$ are the corresponding subalgebras of
strictly triangular matrices and ${{\frak {h}}}$ is the subalgebra of
diagonal matrices. The set ${{\frak {N}}}_{+}\subset {{\frak {gl}}}_{q}$ of
matrices of the form $T=1+S,$ $S\in {{\frak {n}}}_{+},$ is an
infinite-dimensional Lie group with Lie algebra ${{\frak {n}}}_{+}.$

Let us fix  $\lambda \in {{\Bbb {C}}}$ and consider the
evaluation map: 
\begin{equation}
i_{\lambda }:{{\frak {gl}}}_{\infty }\left( {\cal A}_{0}\right) \rightarrow {%
{\frak {gl}}}_{\infty },\quad A\mapsto A(\lambda ).  \label{e24}
\end{equation}

We shall use the following notation: for a subset $K$ of ${{\frak {gl}}}_{q}$
we denote by $K^{\lambda }$ its image under the evaluation map (\ref{e24}).

Our main object, the algebra ${{\frak {gl}}}_{q}^{\lambda }\subset {{\frak {%
gl}}}_{\infty }$ is the image of the whole ${{\frak {gl}}}_{q}.$ We define
also its subalgebras ${{\frak {b}}}_{\pm }^{\lambda },{{\frak {n}}}_{\pm
}^{\lambda },{{\frak {h}}}^{\lambda }$ and the group ${{\frak {N}}}%
_{+}^{\lambda }.$

The algebra ${{\frak {gl}}}_{q}^{\lambda }$ is ${{\Bbb {Z}}}$-graded with
respect to the set of its subspaces ${\cal S}_{n}^{\lambda },\;n\in {{\Bbb {Z%
}}}.$

For any $n\in {{\Bbb {N}}}$ ${{\frak {gl}}}_{n}$ is naturally embedded into $%
{{\frak {gl}}}_{q}^{\lambda }$ as its left upper block: 
\begin{equation}
\left( 
\begin{array}{cc}
{{\frak {gl}}}_{n} & 0 \\ 
0 & 0
\end{array}
\right) .  \label{e25}
\end{equation}

For a matrix $A\in {{\frak {gl}}}_{q}^{\lambda }$ we denote this upper block
by $A_{\mid n}$.

The ${\cal A}_{0}$-valued trace functional on ${{\frak {gl}}}_{q}$ induces
the ordinary ${{\Bbb {C}}}$-valued trace on ${{\frak {gl}}}_{q}^{\lambda }$:
we must put $t=\lambda $ in (\ref{e23}); it will be denoted by the same
symbol. The restriction of $\mathop{Tr}$ to ${{\frak {gl}}}_{n}\subset {%
{\frak {gl}}}_{q}^{\lambda }$ coincides with the standard matrix trace. The
corresponding invariant inner product on ${{\frak {gl}}}_{q}^{\lambda }$ is
non-degenerate: indeed, a matrix which is orthogonal to the all ${{\frak {gl}%
}}_{n},$ $n\in {\Bbb {N}},$ is zero.

\subsection{Loop algebras $L{{\frak {gl}}}_{q},$ $L{{\frak {gl}}}%
_{q}^{\lambda }$}

Due to infinite dimension of ${{\frak {gl}}}_{q},$ ${{\frak {gl}}}%
_{q}^{\lambda },$ an accurate definition of its loop algebras requires some
work. The definitions below have the aim to ensure

(1) the existence of a generalized trace functional and of the corresponding
invariant inner product on $L{{\frak {gl}}}_{q},$ $L{{\frak {gl}}}%
_{q}^{\lambda };$

(2) the possibility to generalize the cross-section theorem \ref{orbits} to
the cases of $L{{\frak {gl}}}_{q},$ $L{{\frak {gl}}}_{q}^{\lambda }.$

Let ${\frak a}$ be an associative algebra; we shall denote by ${\frak a}%
\left( \left( z^{-1}\right) \right) $ the space of formal Laurent series
with coefficients in ${\frak a}.$ For $A\in {\frak a}\left( \left(
z^{-1}\right) \right) $ we denote by $A^{\left[ m\right] }$ its Laurent
coefficient corresponding to $z^{m}.$

Let us consider the algebra ${{\frak {gl}}}_{\infty }\left( {\cal A}%
_{0}((z^{-1}))\right) ,$ i.e. the algebra of ${{\frak {gl}}}_{\infty }$%
-matrices $A$ whose matrix coefficients $A_{ij}(t,z)$ are formal Laurent series
in $z$ with coefficients in ${\cal A}_{0}.$ For any $m\in {\Bbb Z}$ Laurent
coefficients $A_{ij}^{\left[ m\right] }(t)$ form a matrix $A^{\left[
m\right] }=\{A_{ij}^{\left[ m\right] }(t)\}\in {{\frak {gl}}}_{\infty
}\left( {\cal A}_{0}\right) .$ Note that in general $A^{\left[ m\right] }\neq
0$ for all $m\in {\Bbb Z},$ however, for any fixed $i,j$ there exists an
integer $\tilde{N}(A,i,j)$ such that $A_{ij}^{\left[ m\right] }(t)=0$ if $m>%
\tilde{N}(A,i,j).$

\begin{definition}
\label{d22}Loop algebra $L{{\frak {gl}}}_{q}$ consists of matrices $A\in {%
{\frak {gl}}}_{\infty }\left( {\cal A}_{0}((z^{-1}))\right) $ satisfying the
following conditions:

(1) for any $m\in {\Bbb Z}$ $A^{\left[ m\right] }\in {{\frak {gl}}}_{q};$

(2) there exists an integer $N\left( A\right) $ such that ${\rm reg}%
A^{\left[ m\right] }\leq N\left( A\right) ;$

(3) for any $n\in {\Bbb Z}$ there exists an integer $\tilde{N}\left(
A,n\right) $ such that $A_{i,i+n}^{\left[ m\right] }(t)=0$ if $m>\tilde{N}%
(A,n).$
\end{definition}

\vspace{1pt}The notion of regularity degree can be naturally generalized to
the case of the loop algebra $L{{\frak {gl}}}_{q},$ i.e., the regularity
degree ${\rm reg}A$ of a matrix $A\in L{{\frak {gl}}}_{q}$ is the minimal
possible value of $N\left( A\right) .$

We define the ${\cal A}_{0}((z^{-1}))$-valued trace functional on $L{{\frak {%
gl}}}_{q}$ by the same formula as above. In a similar way we may prove that
the trace satisfies $\mathop{Tr}AB=\mathop{Tr}BA.$

The diagonal grading $\left\{ {\cal S}_{i}\right\} ,$ the evaluation map (%
\ref{e24}), the definitions of subalgebras ${{\frak {b}}}_{\pm }^{\lambda },{%
{\frak {n}}}_{\pm }^{\lambda },{{\frak {h}}}^{\lambda }$ and the group ${%
{\frak {N}}}_{+}^{\lambda }$ have their natural counterparts in the case of $%
L{{\frak {gl}}}_{q}.$ For a subset $K\subset {{\frak {gl}}}_{q}\left(
K^{\lambda }\subset {{\frak {gl}}}_{q}^{\lambda }\right) $ we shall denote
by $LK$ $\left( LK^{\lambda }\right) $ the corresponding subset in $L{{\frak 
{gl}}}_{q}$ $\left( L{{\frak {gl}}}_{q}^{\lambda }\right) .$

\section{Gauge orbits of the upper triangular group and the cross-section
theorem.}

In this section we define the gauge action of the upper triangular group and
describe a cross-section of this action. We consider the case of ${{\frak {gl%
}}}_{q},$ the corresponding assertion for ${{\frak {gl}}}_{q}^{\lambda }$
may be obtained by application of the evaluation map.

Let us denote by ${{\Bbb {Y}}}_{q}\subset L{{\frak {gl}}}_{q}$ the affine
subspace of matrices of the form ${\cal L}=\Lambda +A,$ where $A\in L{{\frak {%
b}}}_{+}$ and 
\begin{equation}
\Lambda =\left( 
\begin{array}{cccc}
0 & 0 & 0 & \ldots \\ 
1 & 0 & 0 & \ldots \\ 
0 & 1 & 0 & \ldots \\ 
\vdots & \vdots & \ddots & \ldots
\end{array}
\right) .  \label{e31}
\end{equation}
We define the gauge action of $L{{\frak {N}}}_{+}$ by 
\begin{equation}
{\cal L}\mapsto {\,}^{h}T\cdot {\cal L}\cdot T^{-1},\quad T\in L{{\frak {N}}}%
_{+}.  \label{g1}
\end{equation}
Evidently, the space ${{\Bbb {Y}}}_{q}$ is preserved by this action.

\begin{theorem}
\label{t11} \hspace*{1cm}\newline

\begin{enumerate}
\item  The gauge action of $L{{\frak {N}}}_{+}$ on ${{\Bbb {Y}}}_{q}$ is
free.

\item  The set of companion matrices, i.e. matrices of the form 
\begin{equation}
\widetilde{{\cal L}}=\left( 
\begin{array}{cccc}
u_{1}\left( t,z\right) & u_{2}\left( t,z\right) & u_{3}\left( t,z\right) & 
\ldots \\ 
1 & 0 & 0 & \ldots \\ 
0 & 1 & 0 & \ldots \\ 
\vdots & \vdots & \ddots & \ldots
\end{array}
\right) ,\quad u_{i}\left( t,z\right) \in {\cal A}_{0}{{\Bbb {((}}}z^{-1})),
\label{qfrob}
\end{equation}
is a cross-section of this action.
\end{enumerate}
\end{theorem}

{\em Proof.} Let $T\in L{{\frak {N}}}_{+}$ be an element which converts $%
{\cal L}\in {{\Bbb {Y}}}_{q}$ into a companion matrix $\widetilde{{\cal L}}.$
The statements of the theorem mean that the equation 
\begin{equation}
^{h}T\cdot {\cal L}=\widetilde{{\cal L}}T  \label{e32}
\end{equation}
has a unique solution. Let us write $T,{\cal L},\widetilde{{\cal L}}$ in the
form: 
\begin{equation*}
{\cal L}=\Lambda +\sum\limits_{i\geq 0}{\cal L}^{(i)},\quad \widetilde{{\cal %
L}}=\Lambda +\sum\limits_{i\geq 0}\widetilde{{\cal L}}^{(i)},\quad
T=1+\sum\limits_{j>0}T^{\left( j\right) },
\end{equation*}
where the superscripts $^{(i)}$ denote as above the $i$-th diagonal
component of the corresponding matrices. Substituting this into (\ref{e32})
we obtain the following infinite sequence of equations:

\begin{equation}
\left\{ 
\begin{array}{l}
^{h}T^{\left( 1\right) }\Lambda -\Lambda T^{\left( 1\right) }=-{\cal L}%
^{\left( 0\right) }+{\widetilde{{\cal L}}}^{\left( 0\right) }, \\ 
^{h}T^{\left( i\right) }\Lambda -\Lambda T^{\left( i\right) }=-{\cal L}%
^{\left( i-1\right) }+{\widetilde{{\cal L}}}^{\left( i-1\right)
}+\sum\limits_{j=1}^{i-1}\left( {\widetilde{{\cal L}}}^{\left( i-j-1\right)
}T^{\left( j\right) }-^{h}T^{\left( j\right) }{\cal L}^{\left( i-j-1\right)
}\right) ,\ i\geq 2.\\
\\
\end{array}
\right.  \label{f1}
\end{equation}
The $i$-th equation in this sequence is an equation for $T^{\left( i\right)
} $ and $u_{i}.$ We must prove that:

1) for all $i\in {\Bbb N}$ the corresponding equation has an unique solution 
$T^{\left( i\right) }\in L{\cal S}_{i},$ $u_{i}\in {\cal A}_{0}((z^{-1}));$

2) there exists an integer $N\left( T\right) $ such that ${\rm reg}%
T^{(i)}\leq N(T)$ for all $i\in {\Bbb {N}}.$

The last condition allows to combine all $T^{\left( i\right) }$ into a
single matrix $T\in L{{\frak {N}}}_{+},$ the regularity degree of $T$ being
equal or less than $N\left( T\right) .$ We shall prove not only that $%
N\left( T\right) $ exists, but also that we may put $N\left( T\right) ={\rm %
reg}{\cal L}.$

We shall prove assertion 1), 2) inductively.

Let us rewrite the first equation of (\ref{f1}) as follows: 
\begin{equation}
\left\{ 
\begin{array}{l}
T_{0,1}\left( t,qz\right) =-{\cal L}_{0,0}\left( t,z\right) +u_{1}\left(
t,z\right) , \\ 
T_{n,n+1}\left( t,qz\right) -T_{n-1,n}\left( t,z\right) =-{\cal L}%
_{n,n}\left( t,z\right) ,\quad n\geq 1.
\end{array}
\right.  \label{f2}
\end{equation}
The base of induction is the following

\begin{lemma}
\label{l21} \hspace*{1cm}\newline
\hspace*{1em}1) The equation (\ref{f2}) has a unique solution $T^{\left(
1\right) }\in L{\cal S}_{1},$ $u_{1}\in {\cal A}_{0}((z^{-1})).$

\noindent \hspace*{1em}2) ${\rm reg}T^{\left( 1\right) }\leq {\rm reg}{\cal L%
}.$
\end{lemma}

$\odot $ From (\ref{f2}) it follows that 
\begin{equation}
T_{n,n+1}\left( t,z\right) =-\sum_{i=0}^{n}{\cal L}_{i,i}\left(
t,q^{i-n-1}z\right) +u_{1}\left( t,q^{-n-1}z\right) .  \label{f3}
\end{equation}

$T^{\left( 1\right) }\in L{\cal S}_{1}$ implies that there exists a ${\cal A}%
_{0}$-function $T^{\left( 1\right) }\left( w,t,z\right) $ which interpolates 
$T_{n,n+1}\left( t,z\right) $ for sufficiently large integer $n>N_{1}$: 
\begin{equation}
T_{n,n+1}\left( t,z\right) =T^{\left( 1\right) }\left( n,t,z\right) .
\label{e33}
\end{equation}
By proposition \ref{p12}, 
\begin{equation}
T^{\left( 1\right) }\left( n,n+1,z\right) =0.  \label{e34}
\end{equation}
Substituting (\ref{e33}) and (\ref{e34}) into (\ref{f3}) we find 
\begin{equation}
u_{1}\left( n+1,z\right) =\sum_{i=0}^{n}{\cal L}_{i,i}\left(
n+1,q^{i}z\right) ,\quad \forall n>N_{1}.  \label{f4}
\end{equation}

There exists a unique function $u_{1}\in {\cal A}_{0}\left( \left(
z^{-1}\right) \right) $ satisfying (\ref{f4}). Indeed, developing (\ref{f4})
in powers of $z$ we obtain the following relation for the coefficients $%
u_{1}^{m}\left( t\right) \in {\cal A}_{0}$: 
\begin{equation}
z^{m}:\quad u_{1}^{m}\left( n+1\right) =\sum_{i=0}^{n}{\cal L}%
_{i,i}^{m}\left( n+1\right) \cdot q^{im},\quad \forall n>N_{1}.  \label{f5}
\end{equation}
By definition of regularity degree, for $i\geq {\rm reg}{\cal L}$ all
coefficients ${\cal L}_{i,i}^{m}\left( t\right) $ can be interpolated by $%
{\cal A}_{0}$-functions of two variables $i,t.$ Therefore, by proposition 
\ref{p11p}, the whole sum in the r.h.s. of (\ref{f5}) also may be
interpolated by a (unique) ${\cal A}_{0}$-function for $\forall n>{\rm reg}%
{\cal L},$ and we may put $N_{1}={\rm reg}{\cal L}$.

Once $u_{1}$ is known, $T^{\left( 1\right) }$ is uniquely defined by the
equation (\ref{f3}) ; we can verify in the same way as above that ${\rm reg}%
T^{\left( 1\right) }\leq {\rm reg}{\cal L}.$ $\blacksquare $

Assume now that the first $l$ equations in (\ref{f1}) have unique solutions $%
T^{\left( i\right) }\in L{\cal S}_{i},$ $u_{i}\in {\cal A}_{0}\left( \left(
z^{-1}\right) \right) ,$ $i=1,\ldots ,l,$ and that ${\rm reg}T^{\left(
i\right) }\leq {\rm reg}{\cal L}.$

The $\left( l+1\right) $-th equation has the form: 
\begin{equation}
T_{n,n+l+1}\left( t,qz\right) -T_{n-1,n+l}\left( t,z\right) =\delta
_{n0}u_{l+1}-F_{n,n+l}^{\left( l\right) }\left( t,z\right) ,  \label{f6}
\end{equation}
where $\delta $ is the Kronecker symbol and $F^{\left( l\right) }\in L{\cal S%
}_{l}$ is defined by 
\begin{equation*}
F^{\left( l\right) }={\cal L}^{\left( l\right) }-\sum\limits_{j=1}^{l}\left( 
\widetilde{{\cal L}}^{\left( l-j\right) }T^{\left( j\right) }-^{h}T^{\left(
j\right) }{\cal L}^{\left( l-j\right) }\right) .
\end{equation*}
It is easy to see that ${\rm reg}F^{\left( l\right) }\leq {\rm reg}{\cal L}.$

The condition $T^{\left( l+1\right) }\in L{\cal S}_{l+1}$ imposes $l+1$
restrictions on $T^{\left( l+1\right) }$: 
\begin{equation}
T^{\left( l+1\right) }\left( n,n+k,z\right) =0,\quad k=1,\ldots l+1.
\label{f8}
\end{equation}
In the same way as above the equation (\ref{f8}), corresponding to $k=l+1,$
uniquely defines the coefficient $u_{l+1}$: 
\begin{equation}
u_{l+1}\left( n+l+1,z\right) =\sum_{i=0}^{n}F_{i,i+l}^{\left( l\right)
}\left( n+l+1,q^{i}z\right) ,\quad \forall n>{\rm reg}{\cal L}.  \label{f9}
\end{equation}
Then we find $T^{\left( l+1\right) }$ from 
\begin{equation}
T_{n,n+l+1}\left( t,z\right) =-\sum_{i=0}^{n}F_{i,i+l}^{\left( l\right)
}\left( t,q^{i-n-1}z\right) +u_{l+1}\left( t,q^{-n-1}z\right) .  \label{f10}
\end{equation}
${\rm reg}F^{\left( l\right) }\leq {\rm reg}{\cal L}$ implies that $%
T_{n,n+l+1}$ can be interpolated by a ${\cal A}_{0}$-function for $n>{\rm reg%
}{\cal L}.$

It remains to verify that the conditions (\ref{f8}) for $1\leq k\leq l$ are
also satisfied. Fix some $k.$ Note that $T^{\left( l+1\right) }\left(
n,n+l+1,z\right) =0$ for any $n>{\rm reg}{\cal L},$ therefore, from (\ref{f6}%
) we find: 
\begin{equation*}
T^{\left( l+1\right) }\left( n+l+1-k,n+l+1,q^{l+1-k}z\right)
=-\sum_{i=1}^{l+1-k}F^{\left( l\right) }\left( n+i,n+l+1,q^{i-1}z\right) .
\end{equation*}
All terms in the r.h.s. are zero; indeed, by construction, $F^{\left(
l\right) }\in L{\cal S}_{l}$ and hence $F^{\left( l\right) }(w,w+j,z)=0$ for
all $w\in {\Bbb C}$ and for all  $j$ satisfying $1\leq j\leq l.$ Put $%
w=n+i,j+l+1-i;$  clearly,  $j$ lies in the prescribed range.  So, $T^{\left(
l+1\right) }\left( n,n+k,z\right) =0$ for any $n>{\rm reg}{\cal L}+l+1.$ The
interpolation property gives $T^{\left( l+1\right) }\left( w,w+k,z\right) =0$
for any $w\in {{\Bbb {C}}}.$ But $T_{n,n+l+1}\left( w,z\right) =T^{\left(
l+1\right) }\left( n,w,z\right) ,$ for all $n>{\rm reg}{\cal L},\ w\in {%
{\Bbb {C}}}$ and hence $T_{n,n+l+1}\left( n+k,z\right) =0$ for any $n>{\rm %
reg}{\cal L},$ as desired. $\blacksquare $

\section{The choice of r-matrix.}

Let us fix $\lambda \in {{\Bbb {C}}}.$ Like in the finite-dimensional case
in order to define the generalized DS-reduction we need to find a Poisson
bracket on $L{{\frak {gl}}}_{q}^{\lambda }$ satisfying the following
conditions:

1) the gauge action of $L{{\frak {N}}}_{+}^{\lambda }$ is admissible;

2) the Poisson bracket of any $L{{\frak {N}}}_{+}^{\lambda }$ -invariant
function $\widehat{\psi },$ $\widehat{\psi }_{\mid {{\Bbb {Y}}}_{q}^{\lambda
}}={\rm const},$ with arbitrary $L{{\frak {N}}}_{+}^{\lambda }$-invariant
function vanishes on ${{\Bbb {Y}}}_{q}^{\lambda }.$

We shall use the notation similar to the one of the section 2: ${\cal P}_{+},%
{\cal P}_{-},{\cal P}_{0}$ are the projection operators onto $L{{\frak {n}}}%
_{+}^{\lambda },$ $L{{\frak {n}}}_{-}^{\lambda },$ $L{{\frak {h}}}^{\lambda
},$ respectively; $r=\frac{1}{2}\left( {\cal P}_{+}-{\cal P}_{-}\right)
+r_{0}{\cal P}_{0},$ $r_{0}\in {\rm End}\left( L{{\frak {h}}}^{\lambda
}\right) ;$ $r_{\pm }=r\pm \frac{1}{2}.$ The invariant product on $L{{\frak {%
gl}}}_{q}^{\lambda }$ is defined by

\begin{equation*}
\left\langle A\left( z\right) ,B\left( z\right) \right\rangle =\int \frac{dz%
}{z}\mathop{Tr}A\left( z\right) B\left( z\right) .
\end{equation*}
As in the finite-dimensional case (see \cite{Sem} and the discussion in
section 2), it may be shown that the Poisson bracket of the form 
\begin{equation}
\left\{ \widehat{\varphi },\widehat{\psi }\right\} =\left\langle
\left\langle \left( 
\begin{array}{cc}
r & -\hat{h}r_{+} \\ 
\hat{h}^{-1}r_{-} & -r
\end{array}
\right) \left( 
\begin{array}{c}
\nabla \widehat{\varphi } \\ 
\nabla ^{\prime }\widehat{\varphi }
\end{array}
\right) ,\left( 
\begin{array}{c}
\nabla \widehat{\psi } \\ 
\nabla ^{\prime }\widehat{\psi }
\end{array}
\right) \right\rangle \right\rangle  \label{f31}
\end{equation}
is invariant with respect to the gauge action, and that, moreover,
the gauge action of $L{{\frak {N}}}_{+}^{\lambda }$ is admissible. We put $%
Z_{\widehat{\varphi }}=$\thinspace $^{\hat{h}^{-1}}\nabla \widehat{\varphi }%
-\nabla ^{\prime }\widehat{\varphi },$ $\bar{Z}_{\widehat{\varphi }}=$%
\thinspace $^{\hat{h}^{-1}}\nabla \widehat{\varphi }+\nabla ^{\prime }%
\widehat{\varphi }$ and rewrite (\ref{f31}) as follows: 
\begin{equation}
\left\{ \widehat{\varphi },\widehat{\psi }\right\} =\left\langle Z_{\widehat{%
\varphi }},\tfrac{1}{2}\bar{Z}_{\widehat{\psi }}-rZ_{\widehat{\psi }%
}\right\rangle .  \label{f32}
\end{equation}
Let us define $\hat{s}\in {\rm End}L{{\frak {h}}}^{\lambda }$ by 
\begin{equation}
\hat{s}{\rm diag}\left( f_{0},f_{1},\ldots \right) ={\rm diag}\left(
f_{1},f_{2},\ldots \right) .  \label{f35}
\end{equation}

\begin{proposition}
\label{pp41}

Condition 2) above is equivalent to the following relation for $r_{0}$: 
\begin{equation}
r_{0}\left( 1-\hat{h}\,\hat{s}\right) f=
\tfrac{1}{2}\left( 1+\hat{h}\,\hat{s}\right) f+\alpha \left( f\right),
\quad \forall f\in \Lambda \left( L{{\cal {S}}}_{1}^{\lambda }
\right) ,  \label{f36}
\end{equation}
where $\alpha \left( \cdot \right) :{{\frak {h}}}^{\lambda }\rightarrow {%
{\Bbb {C}}}\cdot 1\subset {{\frak {h}}}^{\lambda }$ is a linear operator.
\end{proposition}

{\em Proof.} Let $\widehat{\varphi },\widehat{\psi }$ be $L{{\frak {N}}}%
_{+}^{\lambda }$ -invariant functions, $\widehat{\psi }_{\mid {{\Bbb {Y}}}%
_{q}^{\lambda }}={\rm const};$ it is easy to see that $Z_{\widehat{\varphi }%
},Z_{\widehat{\psi }}\in L{{\frak {b}}}_{+}^{\lambda }$ and for any ${\cal L}%
\in {{\Bbb {Y}}}_{q}^{\lambda }$ the gradient $d\widehat{\psi }\left( {\cal L%
}\right) \in L{{\frak {n}}}_{+}^{\lambda },$ which implies that $\nabla 
\widehat{\varphi },\nabla ^{\prime }\widehat{\varphi },\bar{Z}_{\widehat{%
\psi }}\in L{{\frak {b}}}_{+}^{\lambda }$ on ${{\Bbb {Y}}}_{q}^{\lambda }.$
Taking into account that $L{{\frak {n}}}_{+}^{\lambda }$ is isotropic and $%
\left( L{{\frak {h}}}^{\lambda }\right) ^{\ast }\simeq L{{\frak {h}}}%
^{\lambda },$ we obtain 
\begin{equation*}
\left\{ \widehat{\varphi },\widehat{\psi }\right\} \left( {\cal L}\right)
=\left\langle Z_{\widehat{\varphi }}^{0},\tfrac{1}{2}\bar{Z}_{\widehat{\psi }%
}^{0}-r_{0}Z_{\widehat{\psi }}^{0}\right\rangle ,\quad {\cal L}\in {{\Bbb {Y}%
}}_{q}^{\lambda },
\end{equation*}
where $Z_{\widehat{\varphi }}^{0}={\cal P}_{0}Z_{\widehat{\varphi }}.$ But $%
\left\{ \widehat{\varphi },\widehat{\psi }\right\} =0$ on ${{\Bbb {Y}}}%
_{q}^{\lambda },$ hence $\tfrac{1}{2}\bar{Z}_{\widehat{\psi }}^{0}-r_{0}Z_{%
\widehat{\psi }}^{0}$ is orthogonal to all $Z_{\widehat{\varphi }}^{0}.$

\begin{lemma}
\label{ll41} For any $f\in L{{\frak {h}}}^{\lambda }$ with $\int \frac{dz}{z}%
\mathop{Tr}f=0$ there exists an $L{{\frak {N}}}_{+}^{\lambda }$-invariant
function $\widehat{\varphi }$ and ${\cal L}\in {{\Bbb {Y}}}_{q}^{\lambda }$
such that $Z_{\widehat{\varphi }}^{0}=f.$
\end{lemma}

This lemma implies that 
\begin{equation}
r_{0}Z_{\widehat{\psi }}^{0}=
\tfrac{1}{2}\bar{Z}_{\widehat{\psi }}^{0}+\tilde{\alpha}\left( d
\widehat{\psi }\right),  \label{f33}
\end{equation}
where $\tilde{\alpha}\left( \cdot \right) :L{{\frak {gl}}}_{q}^{\lambda
}\rightarrow {{\Bbb {C}}}\cdot 1\subset L{{\frak {h}}}^{\lambda }$ is some
linear operator.

Then, it is easy to see that only $L{\cal S}_{1}^{\lambda }$-component of $d%
\widehat{\psi }$ gives contribution in $\bar{Z}_{\widehat{\psi }}^{0},$ $Z_{%
\widehat{\psi }}^{0}$: 
\begin{equation}
\begin{array}{l}
Z_{\widehat{\psi }}^{0}=\left( 1-\hat{h}\,\hat{s}\right) \Lambda d\widehat{%
\psi }^{\left( 1\right) }, \\ 
\bar{Z}_{\widehat{\psi }}^{0}=\left( 1+\hat{h}\,\hat{s}\right) \Lambda d%
\widehat{\psi }^{\left( 1\right) }.
\end{array}
\label{f34}
\end{equation}
Also it is evident that 
\begin{equation}
\tilde{\alpha}\left( d\widehat{\psi }\right) =\alpha \left( \Lambda d%
\widehat{\psi }^{\left( 1\right) }\right) ,  \label{e42}
\end{equation}
where $\alpha \left( \cdot \right) $ is a linear operator in $L{{\frak {h}}}%
^{\lambda }$ with ${\rm Im}\alpha \subset {{\Bbb {C}}}\cdot 1\subset L{%
{\frak {h}}}^{\lambda }.$ Substituting (\ref{f34}) and (\ref{e42}) into (\ref
{f33}) we obtain (\ref{f36}).$\blacksquare $

\begin{remark}
\label{rr41} Proposition \ref{pp41} shows the following important difference
between the ${{{\frak {gl}}}}_{q}^{\lambda }$-case and the
finite-dimensional cases considered in \cite{FRS, SemSev, PirSem}.

In the ${{{\frak {gl}}}}_{n}$-case the diagonal component $\hat{r}^{0}$ of
the r-matrix is given (up to a skew-symmetric operator in $U_{n}$) by the
Cayley transformation of $\hat{h}\tau _{n},$ where $\tau _{n}$ acts in the
subspace of diagonal matrices by cyclic permutation of matrix elements (see (%
\ref{e125})). Obviously, it is impossible to define an analog of $\tau _{n}$
in the ${{{\frak {gl}}}}_{q}^{\lambda }$-case. It is replaced now by the
shift operator $\hat{s},$ whose properties are quite different. This causes
some difficulties. As we shall see below, the operator $A\equiv 1-\hat{h}\,%
\hat{s}$ is not invertible, its kernel is isomorphic to ${{\Bbb {C}}}\left(
\left( z^{-1}\right) \right) $. However, $ImA$ $=L{{{\frak {h}}}}^{\lambda }$
(this is possible only in the infinite-dimensional case), so we can define a
regularized operator $A^{-1}$ and find a r-matrix satisfying (\ref{f36}).
\end{remark}

\vspace{1pt}We shall now study the properties of $A$. We define the
following subspaces in $L{{\frak {h}}}^{\lambda }:$%
\begin{equation}
\begin{array}{l}
V^{\lambda }=\left\{ f\in L{{\frak {h}}}^{\lambda }:f={\rm diag}\left(
0,*,\ldots \right) \right\} ; \\ 
H_{1}^{\lambda }=\left\{ f\in L{{\frak {h}}}^{\lambda }:f\left( \lambda
\right) =0\right\} ; \\ 
V_{1}^{\lambda }=V^{\lambda }\cap H_{1}^{\lambda }; \\ 
{\rm Im}V_{1}^{\lambda }=\left( 1-\hat{h}\,\hat{s}\right) V_{1}^{\lambda };
\\ 
U={\rm Ker}\left( 1-\hat{h}\,\hat{s}\right) .
\end{array}
\label{f37}
\end{equation}

Below we assume that $\lambda \neq 0.$

\begin{proposition}
\label{p31} \hspace*{1cm}\newline
\vspace*{-0.5cm}

\begin{enumerate}
\item  $U=\left\{ {\rm diag}\left( F_{0}\left( z\right) ,F_{0}\left(
q^{-1}z\right) ,\ldots \right) ,\quad F_{0}\left( z\right) \in {{\Bbb {C}}}%
\left( \left( z^{-1}\right) \right) \right\} .$

\item  The restriction of the operator $A$ to $V^{\lambda }$%
\begin{equation*}
A_{\mid V^{\lambda }}:V^{\lambda }\rightarrow L{{\frak {h}}}^{\lambda }
\end{equation*}
is a bijection. We shall denote $A^{-1}\equiv \left[ A_{\mid V^{\lambda
}}\right] ^{-1};$ it is given by 
\begin{equation}
\left( A^{-1}F\right) _{n}\left( z\right) =-\sum_{i=0}^{n-1}F_{i}\left(
q^{i-n}z\right) ,\quad \forall F={\rm diag}\left( F_{0}\left( z\right)
,F_{1}\left( z\right) ,\ldots \right) ,\quad n\in {{\Bbb {N}}}.  \label{f310}
\end{equation}

\item  $V^{\lambda }=V_{1}^{\lambda }\dotplus A^{-1}U.$

\item  With respect to the invariant inner product $L{{\frak {h}}}^{\lambda }
$ is the orthogonal sum 
\begin{equation}
L{{\frak {h}}}^{\lambda }={\rm Im}V_{1}^{\lambda }\oplus U;  \label{f311}
\end{equation}
the projection operator on $U$ is given by 
\begin{equation}
\left( {\cal P}_{U}F\right) _{0}\left( z\right) =-\frac{1}{\lambda }\left(
A^{-1}F\right) \left( \lambda ,q^{\lambda }z\right) ,  \label{f39}
\end{equation}
(obviously, an element of $U$ is uniquely defined by its 0-th component).
\end{enumerate}
\end{proposition}

{\em Proof.} Assertions 1, 2 are evident, they result directly from the
definition of $A.$ To prove~3 suppose that there exists $f\in V_{1}^{\lambda
}\cap A^{-1}U.$ Put $f=A^{-1}F;$ obviously, $F\in U.$ Assertion 1 and (\ref
{f310}) imply that $\left( A^{-1}F\right) _{n}\left( z\right) =-nF_{0}\left(
q^{-n}z\right) $ and hence 
\begin{equation}
\left( A^{-1}F\right) \left( \lambda ,z\right) =-\lambda F_{0}\left(
q^{-\lambda }z\right) .  \label{f312}
\end{equation}
But $\left( A^{-1}F\right) \equiv f\left( \lambda ,z\right) =0,$ since 
$f\in V_{1}^{\lambda },$ hence $F_{0}=0,$ and $f=0.$ So, the sum $%
V_{1}^{\lambda }\dotplus A^{-1}U$ is direct. Let us prove that $V^{\lambda
}\subset V_{1}^{\lambda }\dotplus A^{-1}U$ (the inclusion $\supset $ is
evident). For $f\in V^{\lambda }$ we choose $F\in U$ defined by its
component $F_{0}\left( z\right) =-\frac{1}{\lambda }f\left( \lambda
,q^{\lambda }z\right) .$ The (\ref{f312}) implies that $\left(
f-A^{-1}F\right) \left( \lambda ,z\right) =0,$ i.e. $f-A^{-1}F\in
V_{1}^{\lambda }$ as desired.

Assertions 2 and 3 imply that the sum (\ref{f311}) is direct, formula (\ref
{f39}) follows directly from (\ref{f312}). It remains to verify that ${\rm Im%
}V_{1}^{\lambda }\bot U.$ Recall the definition of $\mathop{Tr}$ on $L{%
{\frak {h}}}^{\lambda }$: $\mathop{Tr}f\left( z\right) ={{\Bbb {D}}}%
_{f}\left( \lambda ,z\right) ,$ where ${{\Bbb {D}}}_{f}\left( \lambda
,z\right) $ is an ${\cal A}_{0}$-function uniquely defined by ${{\Bbb {D}}}%
_{f}\left( n,z\right) =\sum_{i=0}^{n-1}f_{i}\left( z\right) .$ Let $f\in
V_{1}^{\lambda },$ $g\in U.$ We have 
\begin{eqnarray*}
\left\langle Af,g\right\rangle =\left\langle \left( 1-\hat{h}\,\hat{s}%
\right) f,g\right\rangle &=&\int \frac{dz}{z}\left[ \mathop{Tr}fg-\mathop{Tr}%
\left( \hat{s}f\cdot \;^{h^{-1}}g\right) \right] \\
&=&\int \frac{dz}{z}\left[ {{\Bbb {D}}}_{fg}-{{\Bbb {D}}}_{\hat{s}f\cdot
\;^{h^{-1}}g}\right] \left( \lambda ,z\right) .
\end{eqnarray*}
But 
\begin{eqnarray*}
\left[ {{\Bbb {D}}}_{fg}-{{\Bbb {D}}}_{\hat{s}f\cdot \;^{h^{-1}}g}\right]
\left( n,z\right) &=&\sum_{i=0}^{n-1}\left[ f_{i}\left( z\right) g_{i}\left(
z\right) -f_{i+1}\left( z\right) g_{i}\left( q^{-1}z\right) \right] \\
&=&\sum_{i=0}^{n-1}f_{i}\left( z\right) \left[ g_{i}\left( z\right)
-g_{i-1}\left( q^{-1}z\right) \right] \quad \text{(because }f_{0}=0\text{)}
\\
&=&\sum_{i=0}^{n-1}f_{i}\left( z\right) \left[ g_{0}\left( q^{-i}z\right)
-g_{0}\left( q^{-\left( i-1\right) }q^{-1}z\right) \right] =0\quad
\end{eqnarray*}
as desired. $\blacksquare $

\vspace{1pt}

\begin{remark}
\label{rr41p} As above (see remark \ref{rr13}), it can be proved that the
terms containing $\alpha $ do not affect  the corresponding Poisson bracket.
Below  we put $\alpha =0.$
\end{remark}

Taking this into account, we have the following

\begin{theorem}
\label{t33}

A Poisson bracket on $Lgl_{q}^{\lambda }$ of the form (\ref{f31}) admits the
generalized DS-reduction if and only if the corresponding r-matrix $r_{0}$ is
chosen in the form 
\begin{equation}
r_{0}=-\tfrac{1}{2}+A^{-1}+\left( \bar{B}+\tfrac{\lambda }{2}\right) {\cal P}%
_{U},  \label{f325}
\end{equation}
where $\bar{B}$ is a skew-symmetric linear operator in $U.$
\end{theorem}

\begin{remark}
\label{r31}

The subspace $U$ can be naturally identified with ${{\Bbb {C}}}\left( \left(
z^{-1}\right) \right) $ (by taking the 0-th component $F_{0}\left( z\right) $
of  $F\in U$). Hence we may consider the operator $\bar{B}$ as a
skew-symmetric linear operator in ${{\Bbb {C}}}\left( \left( z^{-1}\right)
\right) .$
\end{remark}

{\em Proof of theorem \ref{t33}}. It is easy to see that $\Lambda \left( L%
{\cal S}_{1}^{\lambda }\right) =V_{1}^{\lambda }.$ Then proposition \ref{p31}
shows that any operator $r_{0}$ satisfying (\ref{f36}) can be written in the
form 
\begin{equation}
r_{0}=\frac{1}{2}\left( 1+\hat{h}\,\hat{s}\right) A^{-1}+\tilde{B}{\cal P}%
_{U},  \label{f313}
\end{equation}
where $\tilde{B}$ is a linear operator, $\tilde{B}:U\rightarrow L{{\frak {h}}%
}^{\lambda }.$ The theorem follows directly from

\begin{proposition}
\label{p32} The skew-symmetry of $r_{0}$ is equivalent to the following
conditions:

1) ${\rm Im}\tilde{B}\subset U,$

2) $\tilde{B}+\tilde{B}^{*}=\lambda .$
\end{proposition}

{\em Proof.} The skew-symmetry of $r_{0}$ means that for any $f,g\in L{%
{\frak {h}}}^{\lambda }$%
\begin{equation}
\left\langle r_{0}f,g\right\rangle +\left\langle f,r_{0}g\right\rangle =0.
\label{f315}
\end{equation}
By proposition \ref{p31}, any $f\in L{{\frak {h}}}^{\lambda }$ has a unique
decomposition of the form 
\begin{equation}
f=A\hat{f}+\bar{f},\quad \hat{f}\in V_{1}^{\lambda },\quad \bar{f}\in U.
\label{f316}
\end{equation}
Substituting this in (\ref{f315}) and taking into account that $\tfrac{1}{2}%
\left( 1+\hat{h}\,\hat{s}\right) A^{-1}\equiv -\tfrac{1}{2}+A^{-1},$ we
obtain 
\begin{equation*}
\left\langle r_{0}f,g\right\rangle +\left\langle f,r_{0}g\right\rangle
=I_{1}\left( \hat{f},\hat{g}\right) +I_{2}\left( \hat{f},\bar{g}\right)
+I_{2}\left( \hat{g},\bar{f}\right) +I_{3}\left( \bar{f},\bar{g}\right) ,
\end{equation*}
where 
\begin{eqnarray}
I_{1}\left( \hat{f},\hat{g}\right) &=&-\left\langle A\hat{f},A\hat{g}%
\right\rangle +\left\langle \hat{f},A\hat{g}\right\rangle +\left\langle A%
\hat{f},\hat{g}\right\rangle ,  \label{f318} \\
I_{2}\left( \hat{f},\bar{g}\right) &=&\left\langle \hat{f},\bar{g}%
\right\rangle +\left\langle A\hat{f},A^{-1}\bar{g}\right\rangle
+\left\langle A\hat{f},\tilde{B}\bar{g}\right\rangle ,  \label{f319} \\
I_{3}\left( \bar{f},\bar{g}\right) &=&-\left\langle \bar{f},\bar{g}%
\right\rangle +\left\langle A^{-1}\bar{f},\bar{g}\right\rangle +\left\langle 
\bar{f},A^{-1}\bar{g}\right\rangle +\left\langle \tilde{B}\bar{f},\bar{g}%
\right\rangle +\left\langle \bar{f},\tilde{B}\bar{g}\right\rangle
\label{f320}
\end{eqnarray}
Recalling that $A\equiv 1-\hat{h}\,\hat{s}$ and $\hat{h}^{*}=\hat{h}^{-1}$
we find 
\begin{equation*}
I_{1}\left( \hat{f},\hat{g}\right) =\left\langle \hat{f},\hat{g}%
\right\rangle -\left\langle \hat{s}\hat{f},\hat{s}\hat{g}\right\rangle .
\end{equation*}

\begin{lemma}
\label{l31}

For any $\hat{f},\hat{g}\in V^{\lambda }$ 
\begin{equation*}
\left\langle \hat{f},\hat{g}\right\rangle -\left\langle \hat{s}\hat{f},\hat{s%
}\hat{g}\right\rangle =-\int \frac{dz}{z}\hat{f}\left( \lambda ,z\right) 
\hat{g}\left( \lambda ,z\right) .
\end{equation*}
\end{lemma}

In our case $\hat{f},\hat{g}\in V_{1}^{\lambda }$ i.e. $\hat{f}\left(
\lambda ,z\right) =\hat{g}\left( \lambda ,z\right) =0,$ hence $I_{1}\left( 
\hat{f},\hat{g}\right) =0,$ as desired.

Now, the skew-symmetry of $r_{0}$ implies that both $I_{2}$ and $I_{3}$
are equal to zero.

The first two terms in (\ref{f319}) vanish; indeed, taking into account that 
$\left\langle A\hat{f},\bar{g}\right\rangle =0,$ we have 
\begin{equation*}
\left\langle \hat{f},\bar{g}\right\rangle +\left\langle A\hat{f},A^{-1}\bar{g%
}\right\rangle =\left\langle \hat{f},\bar{g}\right\rangle +\left\langle A%
\hat{f},A^{-1}\bar{g}\right\rangle -\left\langle A\hat{f},\bar{g}%
\right\rangle =I_{1}\left( \hat{f},A^{-1}\bar{g}\right) =0.
\end{equation*}
Hence, $I_{2}=0$ implies $\left\langle A\hat{f},\tilde{B}\bar{g}%
\right\rangle =0,$ i.e. ${\rm Im}\tilde{B}\perp {\rm Im}V_{1}^{\lambda }.$
Then by proposition \ref{p31} ${\rm Im}\tilde{B}\subset U,$ as desired.

Recall now that by (\ref{f312}) $\left( A^{-1}\bar{f}\right) \left( \lambda
,z\right) =-\lambda \bar{f}_{0}\left( q^{-\lambda }z\right) ,$ hence 
\begin{eqnarray*}
0 &=&I_{3}\left( \bar{f},\bar{g}\right) =-\left\langle \bar{f},\bar{g}%
\right\rangle +\left\langle A^{-1}\bar{f},\bar{g}\right\rangle +\left\langle 
\bar{f},A^{-1}\bar{g}\right\rangle +\left\langle \tilde{B}\bar{f},\bar{g}%
\right\rangle +\left\langle \bar{f},\tilde{B}\bar{g}\right\rangle \\
&=&I_{1}\left( A^{-1}\bar{f},A^{-1}\bar{g}\right) +\left\langle \tilde{B}%
\bar{f},\bar{g}\right\rangle +\left\langle \bar{f},\tilde{B}\bar{g}%
\right\rangle \\
&=&-\int \frac{dz}{z}\left( A^{-1}\bar{f}\right) \left( \lambda ,z\right)
\left( A^{-1}\bar{g}\right) \left( \lambda ,z\right) +\left\langle \tilde{B}%
\bar{f},\bar{g}\right\rangle +\left\langle \bar{f},\tilde{B}\bar{g}%
\right\rangle \\
&=&-\lambda ^{2}\int \frac{dz}{z}\bar{f}_{0}\left( q^{-\lambda }z\right) 
\bar{g}_{0}\left( q^{-\lambda }z\right) +\left\langle \tilde{B}\bar{f},\bar{g%
}\right\rangle +\left\langle \bar{f},\tilde{B}\bar{g}\right\rangle \\
&=&-\lambda \left\langle \bar{f},\bar{g}\right\rangle +\left\langle \tilde{B}%
\bar{f},\bar{g}\right\rangle +\left\langle \bar{f},\tilde{B}\bar{g}%
\right\rangle ,
\end{eqnarray*}
i.e. $\tilde{B}+\tilde{B}^{*}=\lambda .$ $\blacksquare $

\section{Explicit formula for the quotient bracket and uniqueness theorem.}

In this section we give an explicit formula for the brackets obtained via
DS-redaction on the quotient ${{\Bbb {Y}}}_{q}^{\lambda }/LN_{+}^{\lambda }$
which may be identified with the set of q-pseudodifference operators of
complex degree $\lambda .$ We shall see that only one of them satisfies
involutivity condition (\ref{e158}). This bracket coincides with the one
constructed in \cite{PirSem}.

In this section we assume that $\lambda $ is generic, i.e. $\lambda \notin 
\frac{2\pi i}{\ln q}{{\Bbb {Q}}}.$

As shown in section 2, the quotient ${{\Bbb {Y}}}_{q}^{\lambda
}/LN_{+}^{\lambda }$ may be identified with the set ${{\Bbb {Y}}}_{0}$ of
companion matrices, i.e. the matrices of the form 
\begin{equation}
{\widetilde{{\cal L}}}=\left( 
\begin{array}{cccc}
-u_{1}\left( z\right) & -u_{2}\left( z\right) & -u_{3}\left( z\right) & 
\ldots \\ 
1 & 0 & 0 & \ldots \\ 
0 & 1 & 0 & \ldots \\ 
\vdots & \vdots & \ddots & \ldots
\end{array}
\right) ,\quad u_{i}\left( z\right) \in {\Bbb C}((z^{-1})).  \label{f41}
\end{equation}
As an affine space ${{\Bbb {Y}}}_{0}$ is isomorphic to the set $\widehat{{%
{\Bbb {G}}}}_{\lambda }$ of q-pseudodifference operators of a complex
degree $\lambda $: 
\begin{equation}
\widehat{{{\Bbb {G}}}}_{\lambda }=\left\{ L=D^{\lambda }+u_{1}\left(
z\right) D^{\lambda -1}+u_{2}\left( z\right) D^{\lambda -2}+\cdots \right\} .
\label{f42}
\end{equation}
Note also that all $\widehat{{{\Bbb {G}}}}_{\lambda }$ are isomorphic (as
affine spaces) to each other and to $\prod\limits_{i\geq 1}{\Bbb C}((
z^{-1}))$; so 
\begin{equation*}
{{\Bbb {Y}}}_{q}^{\lambda }/LN_{+}^{\lambda }\simeq {{\Bbb {Y}}}_{0}\simeq 
\widehat{{{\Bbb {G}}}}_{\lambda }\simeq \prod\limits_{i\geq 1}{\Bbb C}%
((z^{-1})).
\end{equation*}
We fix the following models of the tangent and cotangent spaces of $\widehat{%
{{\Bbb {G}}}}_{\lambda }$: 
\begin{equation}
\begin{array}{l}
T_{L}\widehat{{{\Bbb {G}}}}_{\lambda }=\left\{ X=\bar{X}D^{\lambda }:\bar{X}%
\in J_{-}\subset \Psi D_{q}\right\} , \\ 
T_{L}^{\ast }\widehat{{{\Bbb {G}}}}_{\lambda }=\left\{ f=D^{-\lambda }\bar{f}%
:\bar{f}\in J_{+}\subset \Psi D_{q}\right\} .
\end{array}
\label{f43}
\end{equation}
The canonical pairing between $T_{L}\widehat{{{\Bbb {G}}}}_{\lambda }$ and $%
T_{L}^{\ast }\widehat{{{\Bbb {G}}}}_{\lambda }$ is given by 
\begin{equation}
\left\langle X,f\right\rangle =\mathop{Tr}\nolimits_{\Psi D_{q}} X\cdot f.
\label{f44}
\end{equation}

The space $Fun\left( \widehat{{{\Bbb {G}}}}_{\lambda }\right) $ of the
smooth functional on $\widehat{{{\Bbb {G}}}}_{\lambda }$ is generated by the
Laurent coefficients $u_{i}^{m}$ of the functions $u_{i}\left( z\right) .$
The left and right gradients of a functional $\varphi \in Fun\left( \widehat{%
{{\Bbb {G}}}}_{\lambda }\right) $ are defined by the usual formulas: 
\begin{equation*}
\nabla \varphi \left( L\right) =Ld\varphi ,\quad \nabla ^{\prime }\varphi
=d\varphi L,\quad d\varphi \left( L\right) \in T_{L}^{*}\widehat{{{\Bbb {G}}}%
}_{\lambda },\quad L\in \widehat{{{\Bbb {G}}}}_{\lambda }.
\end{equation*}
It is easy to see that the left and right gradients contain only integer
powers of $D$ and therefore may be considered as elements of $\Psi D_{q}.$

We consider the class of Poisson brackets on $\widehat{{{\Bbb {G}}}}%
_{\lambda } $ of the form 
\begin{equation}
\left\{ \varphi ,\psi \right\} =\left\langle \left\langle \left( 
\begin{array}{cc}
R+aP_{0} & bP_{0} \\ 
cP_{0} & R+dP_{0}
\end{array}
\right) D\varphi ,D\psi \right\rangle \right\rangle ,  \label{f45}
\end{equation}
where 
\begin{equation*}
D\varphi \equiv \left( 
\begin{array}{c}
\nabla \varphi \\ 
\nabla ^{\prime }\varphi
\end{array}
\right)
\end{equation*}
and similarly for $D\psi ,$ $R=\frac{1}{2}\left( P_{+}-P_{-}\right)$ and $%
a,b,c,d$ are linear operators in $J_{0}\simeq {\Bbb C}((z^{-1}))\subset \Psi
D_{q}$ satisfying the skew-symmetry conditions 
\begin{equation*}
a=-a^{\ast },\quad d=-d^{\ast },\quad c^{\ast }=b.
\end{equation*}

\begin{remark}
\label{r41} Note that for a functional $\varphi \in Fun\left( \widehat{{%
{\Bbb {G}}}}_{\lambda }\right) $ its linear gradient $d\varphi $ is defined
up to an arbitrary element of $D^{-\lambda }J_{\left( -\right) }$; in (\ref
{f43}) we have put $D^{-\lambda }d\varphi \in D^{-\lambda }J_{+},$ but it is
a {\em manually} imposed restriction. The bracket (\ref{f45}) is said to be
well-defined if its value does not depend on the $D^{-\lambda }J_{\left(
-\right) }$-components of $d\varphi ,d\psi .$ It is easy to see that the
bracket (\ref{f45}) is well-defined if and only if 
\begin{equation}
a+\tfrac{1}{2}+bD^{-\lambda }=c+\left( \tfrac{1}{2}+d\right) D^{-\lambda
}=\alpha \mathop{Tr}\cdot ,\quad \alpha \in {{\Bbb {C}}}.  \label{f46}
\end{equation}
\end{remark}

\vspace{1pt}

Let ${\cal P}_{00}$ be the projection operator on the one-dimensional
subspace in $U$ generated by the unity matrix, and ${\cal P}_{0}^{ \prime}=%
{\cal P}_{0}-{\cal P}_{00}.$

\begin{theorem}
\label{t41} Let 
\begin{equation}
r_{0}^{\lambda ,\Delta }=-\tfrac{1}{2}+A^{-1}+\left( \bar{B}^{\lambda
,\Delta }+\tfrac{\lambda }{2}\right) {\cal P}_{U},\quad \bar{B}^{\lambda
,\Delta }=\lambda \left( \frac{1}{2}\frac{1+\hat{h}^{\lambda }}{1-\hat{h}%
^{\lambda }}{\cal P}_{0}^{\prime }+\Delta \right) ,\quad   \label{f47}
\end{equation}
where $\Delta $ is a skew-symmetric operator in $U,$ commuting with $\hat{h}.
$ Let $r^{\lambda ,\Delta }=\frac{1}{2}\left( {\cal P}_{+}-{\cal P}%
_{-}\right) +r_{0}^{\lambda ,\Delta }{\cal P}_{0}$.

The Poisson bracket $\left\{ \cdot ,\cdot \right\} _{\Delta }^{\lambda }$ on 
$L{{\frak {gl}}}_{q}^{\lambda }$ defined by 
\begin{equation}
\left\{ \widehat{\varphi },\widehat{\psi }\right\} _{\Delta }^{\lambda
}=\left\langle \left\langle \left( 
\begin{array}{cc}
r^{\lambda ,\Delta } & -\hat{h}r_{+}^{\lambda ,\Delta } \\ 
\hat{h}^{-1}r_{-}^{\lambda ,\Delta } & -r^{\lambda ,\Delta }
\end{array}
\right) \left( 
\begin{array}{c}
\nabla \widehat{\varphi } \\ 
\nabla ^{\prime }\widehat{\varphi }
\end{array}
\right) ,\left( 
\begin{array}{c}
\nabla \widehat{\psi } \\ 
\nabla ^{\prime }\widehat{\psi }
\end{array}
\right) \right\rangle \right\rangle ,  \label{f47p}
\end{equation}
gives rise via DS-reduction to the following bracket on $\widehat{{{\Bbb {G}}%
}}_{\lambda }$: 
\begin{eqnarray}
\left\{ \varphi ,\psi \right\} _{\lambda }^{\Delta } &=&\left\langle
\left\langle \left( 
\begin{array}{cc}
R+\left( \frac{1}{2}\frac{1+\hat{h}^{\lambda }}{1-\hat{h}^{\lambda }}+\Delta
\right) P_{0}^{\prime } & -\left( \frac{1}{1-\hat{h}^{\lambda }}+\Delta
\right) \hat{h}^{\lambda }P_{0}^{\prime } \\ 
\left( \frac{\hat{h}^{\lambda }}{1-\hat{h}^{\lambda }}+\Delta \right) \hat{h}%
^{-\lambda }P_{0}^{\prime } & R-\left( \frac{1}{2}\frac{1+\hat{h}^{\lambda }%
}{1-\hat{h}^{\lambda }}+\Delta \right) P_{0}^{\prime }
\end{array}
\right) D\varphi ,D\psi \right\rangle \right\rangle ;  \notag \\
&&  \label{f48}
\end{eqnarray}
(here we have identified $EndU$ and $End{\Bbb C}((z^{-1})).$ )
\end{theorem}

Like to the finite-dimensional case, the corresponding r-matrix may be
uniquely fixed by imposing, in addition, the involutivity condition:

\begin{theorem}
\label{t42} There exists a unique bracket of the form (\ref{f47p}) on $L%
{\frak {gl}}_{q}^{\lambda }$ which admits DS-reduction and gives rise to a
Poisson bracket on $\widehat{{\Bbb {G}}}_{\lambda }$ satisfying the
involutivity condition 
\begin{equation}
a+b=c+d.  \label{f48p}
\end{equation}
This bracket coincides with $\left\{ \cdot ,\cdot \right\} _{0}^{\lambda }.$
\end{theorem}

\noindent The proof is similar to the one of theorem \ref{tt15} and will be
omitted. \vspace{0.5cm}

The rest of this section is devoted to the proof of theorem \ref{t41}. The
general idea of the proof is similar to the one of \cite{KhM}: both brackets
(\ref{f47p}) and (\ref{f48}) considered as functions of $\lambda $ are
quotients of two ${\cal A}_{0}$-functions; therefore, it is sufficient to
prove the theorem for all sufficiently large integer $\lambda =N.$ But in
this case the DS-reduction on $L{{\frak {gl}}}_{q}^{\lambda }$ amounts to
the q-deformed DS-reduction on $L{{\frak {gl}}}_{N}$ considered in section
2, for which theorem \ref{tt14} gives formula (\ref{f48}).

Let us define filtrations on the spaces $Fun\left( {{\Bbb {Y}}}_{q}^{\lambda
}\right) $ and $Fun\left( \widehat{{{\Bbb {G}}}}_{\lambda }\right) .$ Recall
that we have denoted by $L{\cal S}_{i}\subset L{{\frak {gl}}}_{q}^{\lambda }$
the subspace of matrices which have only $i$-th non-zero diagonal. Let ${%
{\Bbb {V}}}_{n}=\mathop{\textstyle\bigoplus}\limits_{i=0}^{n-1}L{\cal S}%
_{i}, $ let $\Gamma _{n}:{{\Bbb {Y}}}_{q}^{\lambda }\rightarrow {{\Bbb {V}}}%
_{n}$ be the natural projection and $\Omega _{n}=\Gamma _{n}^{\ast }\left(
Fun\left( {{\Bbb {V}}}_{n}\right) \right) .$ Obviously, 
\begin{equation}
\Omega _{1}\subset \Omega _{2}\subset \ldots \subset \Omega _{n}\subset
\ldots \subset Fun\left( {{\Bbb {Y}}}_{q}^{\lambda }\right) ,\quad %
\mathop{\textstyle\bigcup}\limits_{i\geq 1}\Omega _{i}=Fun\left( {{\Bbb {Y}}}%
_{q}^{\lambda }\right) .  \label{f49}
\end{equation}

\begin{lemma}
\label{l41} For any Poisson bracket of the form (\ref{f47p}) 
\begin{equation*}
\left\{ \Omega _{i},\Omega _{j}\right\} _{\Delta }^{\lambda }\subset \Omega
_{i+j}.
\end{equation*}
\end{lemma}

As it was noted $\widehat{{{\Bbb {G}}}}_{\lambda }\simeq \prod\limits_{i\geq
1}{\Bbb C}((z^{-1}))$ as affine spaces. Let ${{\Bbb {M}}}_{n}=\prod%
\limits_{i=1}^{n}{\Bbb C}((z^{-1})),$ let $\gamma _{n}:\prod\limits_{i\geq 1}%
{\Bbb C}((z^{-1}))\rightarrow \prod\limits_{i=1}^{n}{\Bbb C}((z^{-1}))$ be
the natural projection and $W_{n}=\gamma _{n}^{\ast }\left( Fun\left( {{\Bbb 
{M}}}_{n}\right) \right) .$ We have 
\begin{equation}
W_{1}\subset W_{2}\subset \ldots \subset W_{n}\subset \ldots \subset
Fun\left( \widehat{{{\Bbb {G}}}}_{\lambda }\right) ,\quad %
\mathop{\textstyle\bigcup}\limits_{i\geq 1}W_{i}=Fun\left( \widehat{{{\Bbb {G%
}}}}_{\lambda }\right) .  \label{f410}
\end{equation}

\begin{lemma}
\label{l42}For any Poisson bracket of the form (\ref{f48}) 
\begin{equation*}
\left\{ W_{i},W_{j}\right\} _{\lambda }^{\Delta }\subset W_{i+j}.
\end{equation*}
\end{lemma}

The proof of the cross-section theorem \ref{t11} implies that these
filtrations are consistent with the projection $\pi :{{\Bbb {Y}}}%
_{q}^{\lambda }\rightarrow \widehat{{{\Bbb {G}}}}_{\lambda },$ i.e. $\pi
^{\ast }\left( W_{i}\right) \subset \Omega _{i}.$

Let us choose some functions $\varphi \in W_{i},$ $\psi \in W_{j}.$ Let $%
\widehat{\varphi }=\pi ^{*}\varphi ,$ $\widehat{\psi }=\pi ^{*}\psi ;$ they
are defined on ${{\Bbb {Y}}}_{q}^{\lambda }$ and can be continued to $%
LN_{+}^{\lambda }$-invariant functions on the whole $L{{\frak {gl}}}%
_{q}^{\lambda }$ which will be denoted by the same letters. We need to prove
that for any $L\in \widehat{{{\Bbb {G}}}}_{\lambda }$ and some (and then for
any) ${\cal L}\in \pi ^{-1}\left( L\right) $ 
\begin{equation}
\left\{ \widehat{\varphi },\widehat{\psi }\right\} _{\Delta }^{\lambda
}\left( {\cal L}\right) =\left\{ \varphi ,\psi \right\} _{\lambda }^{\Delta
}\left( L\right) .  \label{f411}
\end{equation}
We choose ${\cal L}\in \pi ^{-1}\left( L\right) $ in the companion form (\ref
{f41}). Since $\left\{ \widehat{\varphi },\widehat{\psi }\right\} _{\Delta
}^{\lambda }\in \Omega _{i+j}$ and $\left\{ \varphi ,\psi \right\} _{\lambda
}^{\Delta }\in W_{i+j}$ we may assume that $L\in {{\Bbb {M}}}_{i+j}.$

\begin{lemma}
\vspace{1pt}\label{l44} If (\ref{f411}) holds for all sufficiently large
integer $\lambda ,$ then it holds for all generic $\lambda \in {{\Bbb {C}}}.$
\end{lemma}

{\em Proof. }Since all $\widehat{{{\Bbb {G}}}}_{\lambda }$ are isomorphic to
each other and to $\prod\limits_{i\geq 1}{\Bbb C}((z^{-1})),$ we may consider 
$\left\{ \varphi ,\psi \right\} _{\lambda }^{\Delta }\left( L\right) $ as a
function of variable $\lambda .$ Functionals $\varphi ,\psi $ depend only on
a finite number of Laurent coefficients $u_{i}^{n},$ hence from (\ref{f48})
it follows that $\left\{ \varphi ,\psi \right\} _{\lambda }^{\Delta }\left(
L\right) $ is the quotient of two ${\cal A}_{0}$-functions.

The bracket $\left\{ \widehat{\varphi },\widehat{\psi }\right\} _{\Delta
}^{\lambda }\left( {\cal L}\right) $ may also be considered as a function of 
$\lambda $ because the set ${{\Bbb {Y}}}_{0}$ of companion matrices is
naturally embedded into $L{{\frak {gl}}}_{q}^{\lambda }$ for all $\lambda
\in {{\Bbb {C}}}$. It is easy to see that this bracket is also the quotient
of two ${\cal A}_{0}$-functions. Then the lemma follows directly from the
interpolation property for ${\cal A}_{0}$-functions (see proposition \ref
{p11}). $\blacksquare $

Let us fix some $\lambda =m>i+j.$ The subspace ${{\Bbb {M}}}_{m}\subset 
\widehat{{{\Bbb {G}}}}_{m}$ is Poisson (see \cite{PirSem}). It is naturally
identified with the affine space of normalized $m$-th order q-difference
operators defined in section 2 which was denoted there by the same letter
(see (\ref{e19})). The functions $\varphi ,\psi $ may be considered as
functions on ${{\Bbb {M}}}_{m}.$ As discussed in section 2, the bracket $%
\left\{ \cdot ,\cdot \right\} _{m}^{\Delta }$ on ${{\Bbb {M}}}_{m}$ can be
obtained via the ordinary q-deformed DS-reduction from $L{{\frak {gl}}}_{m}.$
Let $\rho :{{\Bbb {Y}}}_{m}\rightarrow {{\Bbb {Y}}}_{m}/LN_{+}\left(
m\right) \simeq {{\Bbb {M}}}_{m}$ be the corresponding projection, let $\bar{%
\varphi},\bar{\psi}$ be $LN_{+}\left( m\right) $-invariant functions on $L{%
{\frak {gl}}}_{m}$ such that $\bar{\varphi}_{\mid {{\Bbb {Y}}}_{m}}=\rho
^{*}\varphi ,$ $\bar{\psi}_{\mid {{\Bbb {Y}}}_{m}}=\rho ^{*}\psi .$ Theorem 
\ref{tt14} says that 
\begin{equation*}
\left\{ \bar{\varphi},\bar{\psi}\right\} _{\Delta ,m}\left( {{\Bbb {L}}}%
\right) =\left\{ \varphi ,\psi \right\} _{m}^{\Delta }\left( L\right) ,
\end{equation*}
where ${\Bbb L}\in L{{\frak {gl}}}_{m}$ is a companion matrix corresponding
to $L.$ Hence, we need to prove only that 
\begin{equation}
\left\{ \widehat{\varphi },\widehat{\psi }\right\} _{\Delta }^{m}\left( 
{\cal L}\right) =\left\{ \bar{\varphi},\bar{\psi}\right\} _{\Delta ,m}\left( 
{{\Bbb {L}}}\right) .  \label{f412}
\end{equation}

Recall that we have defined an embedding of ${{\frak {gl}}}_{m}$ into ${%
{\frak {gl}}}_{q}^{\lambda }$ (see (\ref{e25})); we may naturally define a
similar embedding for the corresponding loop algebras. As above, for an
element $A\in L{{\frak {gl}}}_{q}^{\lambda }$ we denote by $A_{\mid m}$ its $%
{{\frak {gl}}}_{m}$-block. The similar notation will be used also for
subspaces: for $K\subset L{{\frak {gl}}}_{q}^{\lambda }$ we put 
\begin{equation*}
K_{\mid m}=\left\{ A_{\mid m}:\forall A \in K\right\}.
\end{equation*}

The linear gradient $d\widehat{\varphi }$ is defined up to an arbitrary
element of $L{{\frak {n}}}_{+}^{m};$ however, this freedom does not affect
the value of the bracket. We shall assume that $\left( d\widehat{\varphi }%
\right) _{+}=0,$ so 
\begin{equation}
d\widehat{\varphi }\in \mathop{\textstyle\bigoplus}\limits_{k=-i}^{0}L{\cal S%
}_{k}\subset L{{\frak {b}}}_{-}^{m},  \label{f413p}
\end{equation}
because $\widehat{\varphi }\in \Omega _{i}.$ For $d\bar{\varphi}$ we shall
also use the assumption $\left( d\bar{\varphi}\right) _{+}=0.$

\begin{lemma}
\vspace{1pt}\label{l45} 
\begin{equation}
\left\{ \widehat{\varphi },\widehat{\psi }\right\} _{\Delta }^{m}\left( 
{\cal L}\right) =\left\langle Z_{\bar{\varphi}}\left( {{\Bbb {L}}}\right) ,%
\tfrac{1}{2}\bar{Z}_{\bar{\psi}}\left( {{\Bbb {L}}}\right) -\left(
r^{m,\Delta }\right) _{\mid m}Z_{\bar{\psi}}\left( {{\Bbb {L}}}\right)
\right\rangle .  \label{f415}
\end{equation}
\end{lemma}

{\em Proof. }It is easy to see that 
\begin{equation}
d\widehat{\varphi }_{\mid m}\left( {\cal L}\right) =d\bar{\varphi}\left( {%
{\Bbb {L}}}\right) ;  \label{f413}
\end{equation}
(this follows from the fact that both in ${{\frak {gl}}}_{m}$- and $\ $in ${%
{\frak {gl}}}_{q}^{\lambda }$-cases the gauge action of the upper triangular
group is free and that the restriction of the ${{\frak {gl}}}_{q}^{\lambda } 
$-trace to ${{\frak {gl}}}_{m}$ coincides with the ordinary matrix trace).

Using (\ref{f413p}), (\ref{f413}) and the evident relation ${\cal L}_{\mid
m}={{\Bbb {L}}}$ we find by direct computation that left and right gradients 
$\nabla \widehat{\varphi },\nabla ^{\prime }\widehat{\varphi }$ have the
form 
\begin{equation*}
\left( 
\raisebox{-2.5cm}{
\begin{picture}(50,50)
\put(0,30){\line(1,0){50}} 
\put(20,0){\line(0,1){50}}
\put(17,30){\line(1,-1){30}}
\put(5,30){\line(1,-1){30}}
\put(0,38){\makebox[2cm][r]{$A \in L{\frak{gl}}_m\,$}} \multiput(20,26)(6,-6){4}{\makebox[6mm]{0}} \multiput(45,5)(2,-2){3}{\makebox[2mm]{.}} \put(27.5,9.5){\Large *} \put(15,22){\Large *} \put(10,12){\Huge 0} \put(35,40){\Huge 0} \put(40,20){\LARGE 0} \end{picture}}%
\right) ,
\end{equation*}
and that $\nabla \widehat{\varphi }_{\mid m}=\nabla \bar{\varphi},$ $\nabla
^{\prime }\widehat{\varphi }_{\mid m}=\nabla ^{\prime }\bar{\varphi}.$

Obviously, $Z_{\widehat{\varphi }},$ $\bar{Z}_{\widehat{\varphi }}$ have a
similar form; but $Z_{\widehat{\varphi }}$ is upper-triangular, hence 
\begin{equation}
Z_{\widehat{\varphi }}=\left( 
\begin{array}{ll}
Z_{\bar{\varphi}} & 0 \\ 
0 & 0
\end{array}
\right) .  \label{f414}
\end{equation}
Substituting this into 
\begin{equation*}
\left\{ \widehat{\varphi },\widehat{\psi }\right\} _{\Delta }^{m}\left( 
{\cal L}\right) =\left\langle Z_{\widehat{\varphi }},\tfrac{1}{2}\bar{Z}_{%
\widehat{\psi }}-r^{m,\Delta }Z_{\widehat{\psi }}\right\rangle
\end{equation*}
we obtain (\ref{f415}). $\blacksquare $

\vspace{1pt}

On the other hand, by definition, 
\begin{equation}
\left\{ \bar{\varphi},\bar{\psi}\right\} _{\Delta ,m}\left( {{\Bbb {L}}}%
\right) =\left\langle Z_{\bar{\varphi}},\tfrac{1}{2}\bar{Z}_{\bar{\psi}}-%
\hat{r}_{\Delta ,m}Z_{\bar{\psi}}\right\rangle ,  \label{f416}
\end{equation}
where 
\begin{equation}
\hat{r}_{\Delta ,m}=\frac{1}{2}\left( {\cal P}_{+}-{\cal P}_{-}\right) +\hat{%
r}_{\Delta ,m}^{0},\quad \hat{r}_{\Delta ,m}^{0}=\frac{1}{2}\frac{1+\hat{h}%
\tau _{m}}{1-\hat{h}\tau _{m}}{\cal P}_{0}^{\prime }+m\Delta {\cal P}%
_{U_{m}}.  \label{f416p}
\end{equation}
(Recall that $U_{m}=\left\{ {\rm diag}\left( f_{0}\left( z\right)
,f_{0}\left( q^{-1}z\right) ,\ldots ,f_{0}\left( q^{-\left( m-1\right)
}z\right) \right) \right\} =U_{\mid m}$). The following lemma finishes our
proof:

\begin{lemma}
\label{l43} 
\begin{equation}
\left( r_{0}^{m,\Delta }\right) _{\mid m}=\hat{r}_{\Delta ,m}^{0}.
\label{f416p1}
\end{equation}
\end{lemma}

{\em Proof.} The r-matrix $r_{0}^{m,\Delta }$ satisfies the equation 
\begin{equation}
\tfrac{1}{2}\left( 1+\hat{h}\,\hat{s}\right) f=r_{0}^{m,\Delta }\left( 1-%
\hat{h}\,\hat{s}\right) f,\quad \forall f\in V_{1}^{m}.  \label{f417}
\end{equation}
Evidently, the operator $\hat{h}\,\hat{s}$ preserves $L{{\frak {h}}}_{m},$
which is considered as a subspace of $L{{\frak {h}}}^{m},$ hence (\ref{f417}%
) can be restricted to $L{{\frak {h}}}_{m}$: 
\begin{equation}
\tfrac{1}{2}\left( 1+\hat{h}\,\tau _{m}\right) f=\left( r_{0}^{m,\Delta
}\right) _{\mid m}\left( 1-\hat{h}\,\tau _{m}\right) f,\quad \forall f\in
V_{m},  \label{f418}
\end{equation}
where $V_{m}=\left\{ {\rm diag}\left( 0,*,\ldots *\right) \subset L{{\frak {h%
}}}_{m}\right\} ,\quad V_{m}=\left( V_{1}^{m}\right) _{\mid m}.$ Obviously, $%
\left( r_{0}^{m,\Delta }\right) _{\mid m}$ is skew-symmetric; then,
according to lemma \ref{lll2}, equation (\ref{f418}) implies that 
\begin{equation}
\left( r_{0}^{m,\Delta }\right) _{\mid m}=\frac{1}{2}\frac{1+\hat{h}\tau _{m}%
}{1-\hat{h}\tau _{m}}{\cal P}_{0}^{\prime }+m\tilde{\Delta}{\cal P}%
_{U_{m}}\equiv \hat{r}_{\tilde{\Delta},m}^{0}  \label{f419p}
\end{equation}
with some $\tilde{\Delta}\in {\rm End}U_{m}.$

To prove that $\Delta =\tilde{\Delta}$ let us calculate the bilinear form $%
\left\langle \left( r_{0}^{m,\Delta }\right) _{\mid m}f,g\right\rangle ,$ $%
f,g\in U_{m}.$ For arbitrary $f\in U_{m}$ let us denote by $\hat{f}$ its
image under the natural embedding $L{\frak {gl}}_m \to L{\frak {gl}}^m_q. $

Evidently, 
\begin{equation}
\left\langle \left( r_{0}^{m,\Delta }\right) _{\mid m}f,g \right\rangle_{L%
{\frak {gl}}_m} =\left\langle r_{0}^{m,\Delta } \hat{f},\hat{g}%
\right\rangle_{L{\frak {gl}}^m_q} .  \label{f420}
\end{equation}
By definition, 
\begin{equation}
r_{0}^{m,\Delta }=-\tfrac{1}{2}+A^{-1}+\tfrac{m}{2}{\cal P}_{U}+m\left( 
\frac{1}{2}\frac{1+\hat{h}^{m}}{1-\hat{h}^{m}}{\cal P}_{0}^{\prime }+\Delta
\right) {\cal P}_{U}.  \label{f421}
\end{equation}
Using explicit formula (\ref{f310}) for $A^{-1}$ we find:

\begin{eqnarray}
\left\langle \left( -\tfrac{1}{2}+A^{-1}\right) \hat{f},\hat{g}\right\rangle
&=&-\int \frac{dz}{z}\sum\limits_{l=0}^{m-1}\left( l+\tfrac{1}{2}\right)
f_{0}\left( q^{-l}z\right) g_{0}\left( q^{-l}z\right)  \notag \\
&=&-\frac{m^{2}}{2}\int \frac{dz}{z}f_{0}\left( z\right) g_{0}\left(
z\right) =-\frac{m}{2}\left\langle f,g\right\rangle .  \label{f422}
\end{eqnarray}
Then, from (\ref{f310}) it follows that 
\begin{equation}
\left( {\cal P}_{U}\hat{f}\right) _{\mid m}=f;  \label{f423}
\end{equation}
indeed, 
\begin{equation*}
\left( {\cal P}_{U}\hat{f}\right) _{0}\left( z\right) =\left[ -\frac{1}{%
\lambda }\left( A^{-1}\hat{f}\right) \left( \lambda ,q^{\lambda }z\right)
\right] _{\lambda =m}=\left[ -\frac{1}{\lambda }\left( -m\right) f_{0}\left(
q^{-\lambda }q^{\lambda }z\right) \right] _{\lambda =m}=f_{0}\left( z\right)
.
\end{equation*}
The formulas (\ref{f422}) and (\ref{f423}) imply that 
\begin{equation}
\left\langle \left( r_{0}^{m,\Delta }\right) _{\mid m}f,g\right\rangle
=\left\langle r_{0}^{m,\Delta }\hat{f},\hat{g}\right\rangle =\left\langle
m\left( \frac{1}{2}\frac{1+\hat{h}^{m}}{1-\hat{h}^{m}} {\cal P}%
_{0}^{\prime}+\Delta \right) f,g\right\rangle .  \label{f419p1}
\end{equation}
Recalling that by theorem \ref{tt13} 
\begin{equation}
\left\langle \left( \frac{1}{2}\frac{1+\hat{h}\tau _{m}}{1-\hat{h}\tau _{m}}%
{\cal P}_{0}^{\prime }\right) f,g\right\rangle =\left\langle \frac{m}{2}%
\frac{1+\hat{h}^{m}}{1-\hat{h}^{m}}{\cal P}_{0}^{\prime }f,g\right\rangle
,\quad \forall f,g\in U_{m},  \label{f419}
\end{equation}
and comparing (\ref{f419p}) with (\ref{f419p1}) we find $\tilde{\Delta}%
=\Delta ,$ as desired. $\blacksquare $


\vspace{1cm}

\begin{thebibliography}{99}
\bibitem{a}  {\rm Adler M.} {\it On a trace functional for formal
pseudodifferential operators and the symplectic structure of the Korteweg-de
Vries type equations.} {\rm Inv. Math. {\bf 50}(1979), 219--248.}

\bibitem{ds}  {\rm Drinfeld V.G., Sokolov V.V.} {\it Lie algebras and
equations of Korteweg-de Vries type.} {\rm Sov. Math. Dokl. {\bf 23} (1981),
457--62; J. Sov. Math. {\bf 30} (1985), 1975--2035.}

\bibitem{feig}  {\rm B.L.Feigin.} {\it Lie algebras ${\frak {gl}}(\lambda )$
and cohomology of a Lie algebra of differential operators.} {\rm Russian
Mathematical Surveys {\bf 43}(1988), no. 2, p.169--170.}

\bibitem{ef}  {\rm Feigin B., Frenkel E.} {\it Affine Lie algebras at the
critical level and Gelfand-Dickey algebras.} {\rm Int. J. Math. Phys. {\bf A7%
}, suppl.A1 (1992), 197--215.}

\bibitem{FM}  {\rm L.Freidel, J.M.Maillet , {\it Quadratic algebras
and integrable systems}, Phys. Lett. {\bf B 262,} (1991), p.278.}

\bibitem{Frenkel}  {\rm E.Frenkel.} {\it Deformations of the KDV hierarchy
and related soliton equations.} {\rm Int. Math. Res. Notices {\bf 2}(1996)
55--76; q-alg/9511003.}

\bibitem{FrRes}  {\rm Frenkel E., Reshetikhin N.} {\it Quantum affine
algebras and deformations of the Virasoro algebra and ${\cal W}$-algebras.} 
{\rm Comm. Math. Phys. {\bf 178}(1996), 237--264; q-alg/9505025.}

\bibitem{FRS}  {\rm E.Frenkel, N. Reshetikhin, M.A.Semenov-Tian-Shansky. 
{\it Drinfeld-Sokolov reduction for difference operators and deformations of
W-algebras. I. The case of Virasoro algebra.} Commun. Math. Phys. {\bf 192}
(1998), 605-629.}

\bibitem{gd}  {\rm Gelfand I.M., Dickey L.A.} {\it A family of Hamiltonian
structures related to nonlinear integrable differential equations.} {\rm %
Preprint no. 136, Inst. Appl. Math. USSR Acad. Sci. 1978 (in Russian),
English transl. in: Collected papers of I. M. Gelfand, Vol. 1. Berlin,
Heidelberg, New York : Springer 1987, pp. 625--646}

\bibitem{KhZ}  {\rm B.Khesin and I.Zakharevich.} {\it Poisson-Lie group of
pseudodifferential symbols and fractional KP-KdV hierarchies.} {---{\rm %
C.R.Acad. Sci. Paris, {\bf 316}(1993), Serie I, p.621--626.}}

\bibitem{KhM}  {\rm B.Khesin and F.Malikov.} {\it Universal Drinfeld-Sokolov
Reduction and Matrices of Complex Size.} {\rm Preprint hep-th/9405116,
Commun.Math.Phys. {\bf 175} (1996) 113.}

\bibitem{LP} {\rm Li, L. C., Parmentier, S.} {\it Nonlinear Poisson 
structures and $r$-matrices.}
{\rm Commun. Math. Phys. {\bf 125} (1989), 545--563. }

\bibitem{Oew}  {\rm W.Oewel.} {\it Poisson brackets for integrable lattice
systems.} {\rm \ Algebraic aspects of integrable systems: in memory of Irene
Dorfman. Progress in non-linear differential equations and their
applications, v.26. Ed. A.S.Fokas and I.M.Gelfand, 1997.}

\bibitem{PirSem}  {\rm A.L.Pirozerski and M.A.Semenov-Tian-Shansky.}{\it %
Generalized q-deformed Gelfand-Dickey structures on the group of
q-pseudodifference operators. }{\rm Preprint math.QA/9811025.}

\bibitem{Sem}  {\rm M.A.Semenov-Tian-Shansky.} {\it Dressing transformations
and Poisson-Lie group actions,} {---{\rm Publ. RIMS, Kyoto University {\bf %
21, No.6}(1985), 1203 -- 1221.}}

\bibitem{Monod}  {\rm M.A.Semenov-Tian-Shansky. {\it Monodromy Map and
Classical r-matrices,} Zapiski Nauchn.Semin. POMI, v.200, 1993,
St.Petersburg, preprint hep-th/9402054.}

\bibitem{SemSev}  {\rm M.A.Semenov-Tian-Shansky and A.V.Sevostyanov. {\it %
Drinfeld-Sokolov reduction for difference operators and deformations of
W-algebras. II. General semisimple case.} Commun. Math. Phys. {\bf 192}
(1998), 631.}
\end{thebibliography}
\end{document}